\documentclass[11pt]{article}
\usepackage{amsmath}
\usepackage{amsfonts}
\usepackage{color}
\newtheorem{theo}{Theorem}[section]
\newtheorem{lem}[theo]{Lemma}
\newtheorem{cor}[theo]{Corollary}

\newcommand{\mysection}[1]{\section{#1} \setcounter{equation}{0}}
\newcommand{\proof}{{\sc Proof.} \quad}
\newcommand{\proofc}{{\sc Proof} \ }
\newcommand{\be}{\begin{equation} \label}
\newcommand{\ee}{\end{equation}}
\newcommand{\bea}{\begin{eqnarray}\label}
\newcommand{\eea}{\end{eqnarray}}
\newcommand{\bas}{\begin{eqnarray*}}
\newcommand{\eas}{\end{eqnarray*}}
\newcommand{\bit}{\begin{itemize}}
\newcommand{\eit}{\end{itemize}}
\newcommand{\qed}{\hfill$\Box$ \vskip.2cm}
\newcommand{\nn}{\nonumber}
\newcommand{\R}{\mathbb{R}}
\newcommand{\N}{\mathbb{N}}
\newcommand{\pO}{\partial\Omega}

\newcommand{\eps}{\varepsilon}

\newcommand{\io}{\int_\Omega}
\newcommand{\bom}{\overline{\Omega}}
\newcommand{\abs}{\\[5pt]}
\newcommand{\tm}{T_{\rm max}}
\newcommand{\tmeps}{T_{{\rm max},\eps}}

\newcommand{\ueps}{u_\eps}

\newcommand{\veps}{v_\eps}

\newcommand{\pa}{\partial}

\newcommand{\us}{{u_*}}
\newcommand{\vs}{{v_*}}
\newcommand{\grad}{\nabla}
\renewcommand{\div}{\nabla \cdot}
\begin{document}
\enlargethispage{10mm}
\title{Global existence and stabilization in a diffusive predator-prey model \\
with population flux by attractive transition}
\author{
Frederic Heihoff\footnote{Corresponding author, fheihoff@math.uni-paderborn.de}\\
{\small Institut f\"ur Mathematik, Universit\"at Paderborn,}\\
{\small 33098 Paderborn, Germany} 
\and
Tomomi Yokota\footnote{yokota@rs.kagu.tus.ac.jp}\\
{\small Department of Mathematics, Tokyo University of Science,}\\
{\small Tokyo 162-8601, Japan}
\medskip
}
\date{\small\today}
\maketitle
\begin{abstract}
\noindent 
The diffusive Lotka--Volterra predator-prey model 
\bas
    	\left\{ \begin{array}{rcll}
	u_t &=& \nabla\cdot \left[ d_1\nabla u + \chi v^2 \nabla \Big(\dfrac{u}{v}\Big)\right] 
            +u(m_1-u+av),
	\qquad & x\in\Omega, \ t>0, \\[2mm]
	v_t &=& d_2\Delta v+v(m_2-bu-v),
	\qquad & x\in\Omega, \ t>0, 
	\end{array} \right.
\eas
is considered in a bounded domain 
$\Omega\subset\R^n$, $n \in\{2,3\}$, under Neumann boundary condition,
where $d_1, d_2, m_1, \chi, a, b$ are positive constants and $m_2$ is a real constant. 
The purpose of this paper is to establish global existence and boundedness of 
classical solutions in the case $n=2$ and global existence of weak solutions in the case 
$n=3$ as well as show long-time stabilization. More precisely, we prove that 
the solutions $(u(\cdot,t), v(\cdot,t))$ converge 
to the constant steady state $(u_*, v_*)$
as $t \to \infty$, where $u_*, v_*$ solves $u_*(m_1-u_*+av_*)=v_*(m_2-bu_*-v_*)=0$ with $u_* > 0$ (covering both coexistence as well as prey-extinction cases). 
\abs
\noindent {\bf Key words:} large time behavior; 
predator-prey model; attractive transitional flux\\
 {\bf MSC 2010:} 35B40 (primary); 35K51, 92C17 (secondary)
\end{abstract}
\newpage
\section{Introduction}\label{intro}

Initiated by the modeling of \emph{dictyostelium discoideum} slime mold using a system of coupled partial differential equations in the seminal paper by Keller and Segel in 1970 (\cite{KS1970}) and spurred on by the subsequent successful mathematical verification of its real world aggregation behavior in the model (\cite{W2013}), in the decades since we have seen more and more applications of this same modeling approach to various system of biological agents whose movement is similarly affected by some outside influence. The versatility of this technique is further made evident by it not only generalizing to different dynamics in similar microscopic systems but even to macroscopic systems, such as the modeling of criminal behavior (\cite{SCRIME}) as well as predator-prey models (\cite{OK2018}). Additionally from a mathematical perspective, these kinds of movement dynamics involving an external stimulus can have a significant destabilizing influence on the dynamics of the system as a whole, making even the question of existence of (potentially generalized) global model solutions as well as the nature of long-time behavior challenging and thus answering said questions an often worthwhile endeavor. Given that answers in this regard also help verify whether a proposed model agrees with reality, it is exactly this type of inquiry we will pursue in this paper for the following problem.
\abs
{\bf Problem.}\quad
Motivated by Oeda--Kuto \cite{OK2018}, we consider the diffusive Lotka--Volterra 
predator-prey model with population flux by attractive transition, as given by 
\bas
    	\left\{ \begin{array}{rcll}
	u_t &=& \nabla\cdot \left[ d_1\nabla u + \chi v^2 \nabla \Big(\dfrac{u}{v}\Big)\right] 
            +u(m_1-u+av),
	\qquad & x\in\Omega, \ t>0, \\[2mm]
	v_t &=& d_2\Delta v+v(m_2-bu-v),
	\qquad & x\in\Omega, \ t>0, 
	\end{array} \right.
\eas
which will be precisely studied as the following initial-boundary value problem with  
the equivalent first equation in this paper:
\be{prob}
    	\left\{ \begin{array}{rcll}
	u_t &=& \nabla\cdot ((d_1+\chi v)\nabla u)- \chi \nabla \cdot (u\nabla v) 
            +u(m_1-u+av),
	\qquad & x\in\Omega, \ t>0, \\[2mm]
	v_t &=& d_2\Delta v+v(m_2-bu-v),
	\qquad & x\in\Omega, \ t>0, \\[2mm]
    & & \hspace*{-13mm}\frac{\pa u}{\pa\nu} = \frac{\pa v}{\pa\nu} = 0,
	\qquad & x\in\pO, \ t>0, \\[2mm]
	& & \hspace*{-13mm}
	u(x,0)=u_0(x), \quad
	v(x,0)=v_0(x),
	\qquad & x\in\Omega  
	\end{array} \right.
\ee
in a bounded domain $\Omega\subset\R^n$, $n\in \{2,3\}$, with smooth boundary $\pO$,
where 
\bas
   d_1, d_2, m_1, \chi, a, b >0, \qquad m_2 \in \R
\eas
are constants, $\frac{\pa}{\pa\nu}$ denotes differentiation with respect to
the outward normal of $\pa\Omega$, 
and the initial data $u_0, v_0 \not\equiv 0$ are suitably regular and nonnegative. 
The problem (\ref{prob}) is a diffusive Lotka--Volterra predator-prey model 
in which unknown functions $u=u(x,t)$ and $v=v(x,t)$ represent the population densities 
of the predator and the prey at location $x \in \bom$ and time $t \ge 0$, respectively.  
In this model the drift-diffusion of the predator consists of the population flux 
$-\left[ d_1\nabla u + \chi v^2 \nabla (u/v)\right]$, 
where the first term shows the usual linear diffusion, while 
the nonlinear diffusion described by the second term 
\bas 
 -\chi \nabla\cdot \left[v^2 \nabla \Big(\dfrac{u}{v}\Big)\right]
 = -\chi \nabla\cdot \left[v\nabla u - u \nabla v\right]
\eas
models an ecological situation, in which predators are prone to move towards 
higher concentration of the prey. 
As typical cases of drift-diffusion in biological models, the cross-diffusion 
in the Shigesada--Kawasaki--Teramoto model and the chemotaxis in the Keller--Segel model 
are well-known. As explained by Okubo--Levin \cite{OL2001}, the transition probabilities 
of each predator are determined by conditions at departure, arrival 
and the middle point. 
The so-called taxis $\chi u \nabla v$ is determined by the difference of the conditions 
between departure and arrival and is considered as the middle case, whereas 
$-\chi v^2 \nabla (u/v)$ in (\ref{prob}) shows that the transition probabilities are 
determined by the condition at arrival (\cite{OL2001}) and the latter case 
has not been mathematically studied much except for some recent work on the topic of stationary 
solutions (\cite{KO-preprint}, \cite{OK2018}). 
In particular, there seems to be no work on global existence 
and asymptotic behavior of non-stationary solutions to the problem (\ref{prob}). 
\abs
{\bf Related works.}\quad 
Before stating our results regarding exactly this matter, let us first briefly widen our perspective somewhat to examine the larger context of prior work in this area.
\abs 
In particular for prey-taxis systems with a more standard taxis term, there has been some classical existence (and sometimes boundedness) theory established under the assumption that the taxis term either vanishes for large values of $u$ (\cite{RW_T2010}, \cite{RW_HZ2015}) or the taxis sensitivity is sufficiently small (\cite{RW_W2017}, \cite{RW_WSW2016}). See also \cite{RW_TW2016} for a discussion of classical solution theory in an indirect prey-taxis model. Relaxing the necessary assumptions, there has also been some foray into the construction of weak solutions in similar settings (\cite{ABN2008}, \cite{RW_W2017}).
\abs
Often already established in tandem with the already mentioned existence theory, there has also been some exploration of long-time behavior in various diffusive predator-prey models. More precisely, it has been shown that in many scenarios the solutions tend (potentially exponentially fast) to their constant equilibria (\cite{RW_JW2017}, \cite{RW_TW2016}, \cite{RW_WSW2016}). In \cite{RW_W2017}, the constructed weak solutions are even further shown to become classical as a consequence of stabilization.
\abs
Moving yet slightly further away from the model we are interested in but staying in the realm of predator-prey modeling and analysis, there have been some efforts made to analyze predator-taxis systems, which model prey fleeing from its predators as opposed to the predators chasing the prey (\cite{RW_WWS2018}). When this dynamic is combined with the already discussed prey-taxis, the resulting system seems to be rather challenging as there are to our knowledge thus far only some constructions of weak solutions available (\cite{RW_TW2020} considering the one dimensional case, \cite{RW_F2021} featuring nonlinear diffusion) and classical solutions seem to have only been constructed if both taxis terms are mediated by a regularizing indirection (\cite{RW_W2022}) or the initial data is already very close to the equilibria (\cite{RW_F2020}). There has also been some discussion of models featuring multiple predators attracted by the same prey (\cite{RW_WW2018}).
\abs
{\bf Purpose.}\quad 
The purpose of this paper is 
\begin{enumerate}
\item[1)] to establish global existence and boundedness of classical solutions 
in the case $n=2$ and global existence of weak solutions in the case $n=3$;
\item[2)] to show stabilization or more specifically to show that the solution $(u(\cdot,t), v(\cdot,t))$ converges 
to the constant steady state $(u_*, v_*)$ in some appropriate sense as $t \to \infty$, 
where $u_*, v_*$ are either the positive solutions to $u_*(m_1-u_*+av_*)=v_*(m_2-bu_*-v_*)=0$ in the coexistence case ($m_2 - bm_1 > 0$) or equal to $m_1$ and $0$, respectively, in the prey-extinction case ($m_2 - bm_1 \leq 0$). 
\end{enumerate}

{\bf Main results.}\quad 
The main result for the first purpose concerning global existence reads as follows. 
\begin{theo}\label{theo_global}
  Let $n \in \{2,3\}$, $d_1, d_2, m_1, \chi, a, b >0$, $m_2 \in \R$ and assume that  
$u_0,v_0 \in W^{1,\infty}(\Omega)$ are nonnegative with $u_0, v_0 \not\equiv 0$. In addition, if $n=3$, assume that $u_0 \in L\log L(\Omega)$ and $\sqrt{v_0} \in W^{1,2}(\Omega)$. 
\begin{enumerate}
\item[\rm (i)] If $n=2$, then the problem (\ref{prob}) possesses a unique global 
classical solution $(u,v)$ such that 
  \bas
\left\{  \begin{array}{l}
	u \in C^0(\bom \times [0,\infty)) \cap C^{2,1}(\bom\times (0,\infty)), \\[1mm]
	v \in C^0(\bom \times [0,\infty)) \cap C^{2,1}(\bom\times (0,\infty)) 
	\end{array} \right.
  \eas
and such that $u,v>0$ in $\bom \times (0,\infty)$. Moreover, this solution is 
bounded in the sense that there exists a constant $C>0$ fulfilling 
  \bas
	\|u(\cdot,t)\|_{L^\infty(\Omega)}+ \|v(\cdot,t)\|_{W^{1,\infty}(\Omega)} \le C 
    \qquad \mbox{for all}\ t>0.
  \eas
\item[\rm (ii)] If $n=3$, then there exists at least one global weak solution 
$(u,v)$ of (\ref{prob}) in the sense that 
  \bas
\left\{  \begin{array}{l}
	u \in L^{\frac{4}{3}}_{\rm loc}([0,\infty); W^{1,\frac{4}{3}}(\Omega)), \\[1mm]
	v \in L^4_{\rm loc}([0,\infty); W^{1,4}(\Omega))
	\end{array} \right.
  \eas
and that $u,v \ge 0$ a.e.\ on $\Omega \times (0,\infty)$, 
  \bas
   u \in L^2_{\rm loc}([0,\infty); L^2(\Omega)), \qquad  
   v \in L^\infty(\Omega \times (0,\infty))
  \eas
as well as  the identities
  \bea{}
  -\int_0^\infty \io u \varphi_t - \io u_0\varphi(\cdot,0)
  &=&-\int_0^\infty \io (d_1+\chi v) \nabla u \cdot \nabla\varphi
   +\chi \int_0^\infty \io u \nabla v \cdot \nabla\varphi\nn \\
  & &\hspace{34.5mm} +\int_0^\infty \io u(m_1-u+av)\varphi, \label{weak_solution_u} \\[2mm]
   -\int_0^\infty \io v \varphi_t - \io v_0\varphi(\cdot,0) 
  &=&-d_2\int_0^\infty \io \nabla v \cdot \nabla\varphi
   +\int_0^\infty \io v(m_2-bu-v)\varphi  \label{weak_solution_v} 
  \eea
hold for all $\varphi \in C_0^\infty(\bom\times[0,\infty))$.
\end{enumerate}
\end{theo}
Given that we have now properly formalized the first main result of this paper, we can now transition to stating our second main result concerning the stabilization behavior of the now established solutions to (\ref{prob}).
\begin{theo}\label{theo_stabilization}
Let $n \in \{2,3\}$, $d_1, d_2, m_1, \chi, a, b >0$, $m_2 \in \R$ and assume that  
$u_0,v_0 \in W^{1,\infty}(\Omega)$ are nonnegative with $u_0, v_0 \not\equiv 0$. In addition, if $n=3$, assume that $u_0 \in L\log L(\Omega)$ and $\sqrt{v_0} \in W^{1,2}(\Omega)$. Let
	\be{equilibria}
		\us := \begin{cases}
			\frac{m_1 + a m_2}{ab + 1} &\text{ if } m_2 - bm_1 \geq 0 \\
			m_1 &\text{ if } m_2 - bm_1 < 0
		\end{cases}
		\;\;\;\;\text{ and }\;\;\;\; \vs := \begin{cases}
			\frac{m_2 - b m_1}{ab + 1} &\text{ if } m_2 - bm_1 \geq 0 \\
			0 &\text{ if } m_2 - bm_1 < 0
		\end{cases}.
	\ee
	be the constant steady states. If
	\be{stabilization_initial_data_condition} 
		\log (u_0) \in L^1(\Omega) \;\;\;\; \text{ and } \;\;\;\; 
       \log (v_0) \in L^1(\Omega) \text{ if } v_\star > 0
	\ee
	as well as
	\be{stablization_param_condition}
	\chi^2 < \frac{4 d_1 d_2}{b m_{2+} \us } \left( \frac{a \vs}{m_{2+}} +\frac{4}{b} \right), \;\;\;\;  \text{ where } m_{2 +} := \max(0,m_2),
	\ee
	then the solution $(u,v)$ constructed in Theorem \ref{theo_global} has one of the following stabilization properties depending on the space dimension $n$:
	\begin{enumerate} 
	\item[\rm (i)] If $n = 2$, then
	\be{stabilization_2d_1}
		u(\cdot, t) \rightarrow \us \;\; \text{ in } L^p(\Omega) \;\;\;\; \text{ and }\;\;\;\; v(\cdot, t) \rightarrow \vs \;\; \text{ in } W^{1,p}(\Omega)
	\ee
	for all $p \in [1,\infty)$ as $t \rightarrow \infty$. 
	\item[\rm (ii)] If $n = 3$, then
	\be{stabilization_2d_2}
		u(\cdot, t) \rightarrow \us \;\; \text{ in } L^1(\Omega) \;\;\;\; \text{ and }\;\;\;\; v(\cdot, t) \rightarrow \vs \;\; \text{ in } L^p(\Omega) \;\;\;
	\ee
	for all $p \in [1,\infty)$ as $(0, \infty) \setminus N \ni t \rightarrow \infty$, where $N$ is a set of measure zero. 
	\end{enumerate}
\end{theo}
{\bf Main ideas.}\quad
Before diving into the details of the proofs, we will first give a brief overlook over our main ideas as well as the organization of this paper. 
\abs
In Section \ref{sec: local}, we begin by establishing some basic local existence theory for not only (\ref{prob}), which we will directly investigate for the two-dimensional case, but also for a family of variants featuring a regularization of the taxis term, which we will use in the three-dimensional case to gain our desired weak solution as a limit of the thus gained approximate solutions. We further derive a set of uniform baseline a priori estimates for these local solutions.
\abs
Having established the necessary preliminaries, we devote Section \ref{sec: global} to the proof of Theorem \ref{theo_global}. To this end, we begin by treating the two-dimensional case following the approach from \cite{TW2016_ZAMP} by essentially iteratively bootstrapping ourselves to sufficiently good a priori estimates to rule out finite-time blow-up in (\ref{prob}). To handle the three-dimensional case, we do not work directly with (\ref{prob}) because of the lack of regularity but instead with the family of similar systems seen in (\ref{prob-ap}) only differing due to a slight regularization of the taxis term. Our aim is then to gain our desired solution as the limit of the thus gained family of approximate solutions.
The key to establishing the necessary bounds to achieve this is the derivation of an energy-type inequality of the form 
\bas
  & &\hspace{-6mm} \frac{d}{dt}\left[\io u\log(u) 
          + \frac{\chi}{2b} \io \frac{|\nabla v|^2}{v} \right]
      + \frac{1}{K}\left[
        \io \frac{d_1+\chi v}{u}|\nabla u|^2 
      + \io u^2 \log(u)
      + \io \frac{|D^2 v|^2}{v}
      + \io \frac{|\nabla v|^4}{v^3}\right] \nn
  \le K 
\eas
in Lemma \ref{lem_energy1}. We then use the bounds resulting from the above inequality as well as some of their immediate consequences to both ensure that the approximate solutions are in fact global in Lemma \ref{lem_global-app} by employing similar arguments to the ones used in the two-dimensional case as well as to construct our actual solution candidates by using compactness arguments in Lemma \ref{lem_limit-process}. It then only remains to show that the solution properties of the approximate solutions survive the limit process in such a way as to make our newly constructed solution candidates actual weak solutions to (\ref{prob}).
\abs
Given that we have now firmly established the existence of our desired solutions, Section \ref{sec: stabilize} focuses on the analysis of their long-time behavior and thus the proof of Theorem \ref{theo_stabilization}. Here, we will not only treat both the two-dimensional and three-dimensional cases at once but also deal with both the coexistence and prey-extinction cases in a way that tries to minimize the duplication of similar arguments. We do this by basing our arguments for all cases on another energy-type inequality, which has the form
\bas
\frac{d}{dt} E(u,v) \leq -\delta G(u,v), \quad\text{ where }\quad E(u,v) := \int_\Omega H_\us(u) + \frac{a}{b}\int_\Omega H_\vs(v) + \frac{2}{b^2 \widehat{m_2}}\int_\Omega (v - \vs)^2
\eas 
with $H_\xi(\eta) = \eta - \xi - \xi \log(\frac{\eta}{\xi})$ if $\xi > 0$ and $H_0(\eta) = \eta$, and
\bas
  G(u,v) := \int_\Omega \frac{|\grad u|^2}{u^2} +\int_\Omega \frac{|\grad v|^2}{v^2} + \int_\Omega (u - \us)^2 + \int_\Omega (v - \vs)^2.
\eas
We then show in Lemma \ref{lem_stab_energy} that, under the parameter condition (\ref{stablization_param_condition}), both the classical solution in the two-dimensional case as well as the approximate solutions in the three-dimensional case fulfill said energy-type inequality after some (uniform) waiting-time $T_E$ while the energy $E$ is (uniformly) bounded up to said same time. The need for this waiting time is owed to eliminiating an otherwise necessary initial data condition by using Lemma \ref{lem_eventual_v_bound}. As a consequence of this energy-type inequality, we conclude that the energy $E$ is monotonically decreasing for almost all times $t > T_E$ in Lemma \ref{lem_E_monotone} by first deriving a space-time square integrability bound for $H_\us(u)$ and $H_\vs(v)$ in Lemma \ref{lem_H2_bounds}, which is necessary to ensure that the monotonicity property survives the limit process used in our weak solution construction, using properties of the functions $H_\xi(\eta)$ established in e.g.\ \cite{W2020_ANS}. Combining this same square integrability property with the second notable consequence of the energy-type inequality,
that is, $\int_{T_E}^\infty G(u,v) < \infty$,
we can construct a sequence of times along which $E(u,v)$ converges to zero, which combined with the already established monotonicity almost immediately yields a slightly weaker version of our desired result in Lemma \ref{lem_stabilization} and after some slight refinement Theorem \ref{theo_stabilization} itself.
\mysection{Preliminaries. Local existence and basic estimates}\label{sec: local}
In this section we give lemmas which present local existence and 
basic estimates for  classical solutions of (\ref{prob}) as well as a more regularized version of the same system, which we will later use in our weak solution construction. 
For convenience we write down the first and second equations of the system in equation as follows:
\be{prob-general}
    	\left\{ \begin{array}{rcll}
	u_t &=& \nabla\cdot ((d_1+\chi v)\nabla u)- \chi \nabla \cdot (F(u)\nabla v) 
            +u(m_1-u+av), \qquad & x\in\Omega, \ t>0, \\[2mm]
	v_t &=& d_2\Delta v+v(m_2-bu-v), \qquad & x\in\Omega, \ t>0, \\[2mm]
    & & \hspace*{-13mm}\frac{\pa u}{\pa\nu} = \frac{\pa v}{\pa\nu} = 0,
	\qquad & x\in\pO, \ t>0, \\[2mm]
	& & \hspace*{-13mm}
	u(x,0)=u_0(x), \quad
	v(x,0)=v_0(x),
	\qquad & x\in\Omega,   
	\end{array} \right.
\ee
where $F$ is given by 
\bas
   F(u):=
   \left\{ \begin{array}{ll}
	u \quad & \mbox{if}\quad n=2, \\[2mm]
	\dfrac{u}{1+ \eps u} \quad (\eps>0) \quad & \mbox{if}\quad n=3.
	\end{array} \right.
\eas

The following statement on local existence and uniqueness can be proved 
by using fixed point arguments as in the proof of 
\cite[Lemma 2.1]{TW2011_SIAM} or \cite[Lemma 2.1]{W2019_JDE}.  
\begin{lem}\label{lem_loc}
  Let $n \in \{2,3\}$, $d_1, d_2, m_1, \chi, a, b >0$, $m_2 \in \R$ and assume that  
$u_0,v_0 \in W^{1,\infty}(\Omega)$ are nonnegative and $u_0, v_0 \not\equiv 0$. 
Then there exist $\tm \in (0,\infty]$ and 
a uniquely determined pair $(u,v)$ of functions
  \bas
	& & u \in C^0(\bom \times [0,\tm)) \cap C^{2,1}(\bom\times (0,\tm)) \qquad \mbox{and} \nn\\
	& & v \in C^0(\bom \times [0,\tm)) \cap C^{2,1}(\bom\times (0,\tm)) 
  \eas
  such that $u,v>0$ in $\bom \times (0,\tm)$,   
  such that $(u,v)$ solves (\ref{prob-general}) classically, and such that 
  \be{ext}
        \mbox{if $\tm<\infty$, \quad then} \quad 
	\|u(\cdot,t)\|_{L^\infty(\Omega)}+ \|v(\cdot,t)\|_{W^{1,\infty}(\Omega)} \to \infty 
    \ \mbox{as}\ t\nearrow\tm.
  \ee
\end{lem}
We now establish some basic estimates that will not only provide the foundation for our arguments refuting blow-up of the now established local solutions but will also prove useful in our later discussion of long-time behavior. Importantly, at no point are the constants derived in the following two lemmas dependent on the exact structure of $F$ as the term containing $F$ in the first equation always vanishes as a result of integration by parts and the prescribed Neumann boundary conditions.
\begin{lem}\label{lem_basic-v}
  Under the assumption of Lemma \ref{lem_loc}, there exists a constant $C>0$ such that 
  \be{v-bound}
    \|v(\cdot,t)\|_{L^\infty(\Omega)} \le C \qquad \mbox{for all}\ t \in (0,\tm)   
  \ee
and 
  \be{nabla-v-bound}
	\int_t^{t+\tau} \io |\nabla v|^2 \le C \qquad \mbox{for all } t\in (0,\tm-\tau)
  \ee
with $\tau:=\min\{1,\frac{1}{2}\tm\}$.
\end{lem}
\proof%
Multiplying the second equation in (\ref{prob-general}) by $v^{p-1}$ ($p>1$) 
and integrating it over $\Omega$, we see from the H\"older inequality that
  \bea{ODI-1}
    \frac{1}{p}\dfrac{d}{dt}\io v^p + d_2(p-1)\io v^{p-2}|\nabla v|^2
    &=& m_2 \io v^p - b \io u v^p + \io v^{p+1} \nn\\ 
    &\le& m_2 \io v^p - |\Omega|^{-\frac{1}{p}}\left(\io v^p\right)^{1+\frac{1}{p}}. 
  \eea
By dropping the second term on the left-hand side, an ODE comparison shows that 
  \bas
    \io v^p(\cdot,t) \le \max\left\{\io v_0^p,\ m_{2+}^p|\Omega|\right\}, \quad \mbox{i.e.,}\quad  
    \|v(\cdot,t)\|_{L^p(\Omega)} 
    \le \max\{\|v_0\|_{L^p(\Omega)},\ m_{2+}|\Omega|^\frac{1}{p}\}. 
  \eas
Letting $p \to \infty$ precisely warrants (\ref{v-bound}). 
Moreover, a straightforward integration of (\ref{ODI-1}) with $p=2$ in time 
yields (\ref{nabla-v-bound}).
\qed
\begin{lem}\label{lem_basic-u}
  Under the assumption of Lemma \ref{lem_loc}, there exists a constant $C>0$ such that 
  \be{mass}
	\io u(\cdot,t) \le C \qquad \mbox{for all } t\in (0,\tm)
  \ee
and 
  \be{u-squared-bound}
	\int_t^{t+\tau} \io u^2 \le C \qquad \mbox{for all } t\in (0,\tm-\tau)
  \ee
with $\tau:=\min\{1,\frac{1}{2}\tm\}$.
\end{lem}
\proof%
Integrating the first equation in (\ref{prob-general}) over $\Omega$, 
we immediately infer from (\ref{v-bound}) that 
  \bea{ODI-2}
   \frac{d}{dt}\io u 
   =  \io (m_1+av)u - \io u^2 
   \le (m_1+ac_1) \io u - \io u^2 
  \eea
for some $c_1>0$. Thus the Cauchy--Schwarz inequality gives 
  \bas
   \frac{d}{dt}\io u 
   \le (m_1+a c_1) \io u - |\Omega|^{-1}\left(\io u\right)^2,
  \eas
which by an ODE comparison derives 
  \bas
     \io u (\cdot,t) \le \max\left\{\io u_0,\ (m_1+ac_1)|\Omega|\right\}.
  \eas
This proves (\ref{mass}), and then (\ref{u-squared-bound}) results from a 
straightforward integration of (\ref{ODI-2}) in time.
\qed
\mysection{Global existence}\label{sec: global}
In this section, following the arguments in \cite{TW2016_ZAMP}, \cite{W2012}   
and \cite{W2016}, we establish global existence in (\ref{prob}). 
The proof will be divided into the two cases that $n = 2$ in Section \ref{sec: 2D} 
and that $n = 3$ in Section \ref{sec: 3D}. 
%
%
%
%
%
%
\subsection{Global existence and boundedness of classical solutions in the case $n=2$}\label{sec: 2D}
We will show that $\tm=\infty$ as in the proof of \cite{TW2016_ZAMP} while we will use an energy inequality approach to show the same for the approximate solutions used in Section \ref{sec: 3D} 
to handle the case $n=3$. 
\begin{lem}\label{lem_Lap-v}
  Under the assumption of Lemma \ref{lem_loc}, there exists a constant $C>0$ such that 
  \be{nabla-v-uniform}
    \io |\nabla v(\cdot,t)|^2 \le C \qquad \mbox{for all}\ t \in (0,\tm)   
  \ee
and 
  \be{Lap-v-bound}
	\int_t^{t+\tau} \io |\Delta v|^2 \le C \qquad \mbox{for all } t\in (0,\tm -\tau)
  \ee
with $\tau:=\min\{1,\frac{1}{2}\tm\}$.
\end{lem}
\proof%
Testing the second equation in (\ref{prob-general}) by $-\Delta v$, we see from 
(\ref{v-bound}) and Young's inequality that 
  \bea{ODI-3}
    \frac{1}{2}\dfrac{d}{dt}\io |\nabla v|^2 + d_2 \io |\Delta v|^2 
    &=& m_2 \io v(-\Delta v) + b \io uv\Delta v - 2 \io v |\nabla v|^2 \nn\\
    &\le& c_1 + c_2 \io u^2 + \frac{d_2}{2} \io |\Delta v|^2 
  \eea
for some $c_1,c_2>0$. 
In view of (\ref{nabla-v-bound}) we know that there exists a constant $c_3$ such that
  \bas
    \int_t^{t+\tau} \io |\nabla v|^2 \le c_3 \qquad \mbox{for all } t\in (0,\tm-\tau),
  \eas
and hence, given $t \in (0,\tm)$, 
we can pick $t_0 \in (t-\tau, t) \cap [0,\infty)$ fulfilling
  \bas
    \io |\nabla v(\cdot,t_0)|^2 \le c_4:=\max\left\{\io |\nabla v_0|^2,\ \frac{c_3}{\tau}\right\}.
  \eas
Integrating (\ref{ODI-3}) over $(t_0,t)$ and invoking (\ref{u-squared-bound}), 
we thus infer (\ref{nabla-v-uniform}), whereupon another integration of (\ref{ODI-3}) 
yields (\ref{Lap-v-bound}). 
\qed
\begin{lem}\label{lem_inductive-u}
  Under the assumption of Lemma \ref{lem_loc}, let $p \ge 2$, $L>0$ and 
$\tau:=\min\{1, \frac{1}{2}\tm\}$. Then there exists a constant $C=C(p,L)>0$ such that 
if  
  \bas
    \int_t^{t+\tau}\io u^p \le L \qquad \mbox{for all}\ t \in (0,\tm-\tau),   
  \eas
then  
  \bas 
	\io u^p(\cdot,t) \le C \qquad \mbox{for all } t\in (0,\tm)
  \eas
and 
  \bas
    \int_t^{t+\tau}\io u^{p+1} \le C \qquad \mbox{for all}\ t \in (0,\tm-\tau). 
  \eas
\end{lem}
\proof%
Testing the first equation in (\ref{prob-general}) by $u^{p-1}$, we see from 
integration by parts and (\ref{v-bound}) that 
  \bas
    \frac{1}{p}\frac{d}{dt}\io u^p + (p-1) \io (d_1+\chi v)u^{p-2}|\nabla u|^2  
    &=& (p-1)\chi \io u^{p-1}\nabla u \cdot \nabla v 
        +\io u^p(m_1-u+av) \\
    &\le& (p-1)\chi \io u^{p-1}\nabla u \cdot \nabla v 
        +c_1\io u^p - \io u^{p+1}
  \eas
for some $c_1>0$. This implies that the situation is the same as 
in \cite[Proof of Lemma 3.8, (3.32)]{TW2016_ZAMP}. 
Therefore the same argument yields the lemma. 
\qed 
Recalling the basic estimate for $u^2$ in Lemma \ref{lem_basic-u}, 
we can now obtain the following lemma immediately.
\begin{lem}\label{lem_L^p-u}
  Under the assumption of Lemma \ref{lem_loc}, for all $p \in (1,\infty)$ 
there exists a constant $C(p)>0$ such that 
  \bas 
	\io u^p(\cdot,t) \le C(p) \qquad \mbox{for all } t\in (0,\tm).
  \eas
\end{lem}
\proof%
The lemma results from an induction via Lemma \ref{lem_inductive-u} 
with the starting estimate (\ref{u-squared-bound}). 
\qed
We are now in a position to prove global existence in the case $n=2$. \abs
\proofc of Theorem \ref{theo_global} (i).\quad
Let $(u,v)$ denote the classical solution of \eqref{prob-general}  
in $\bom \times [0,\tm)$ provided by Lemma \ref{lem_loc}. 
According to Lemmas \ref{lem_basic-v} and \ref{lem_L^p-u}, 
we know that  
  \bas
    \|v(\cdot,t)\|_{L^\infty} \le c_1 \qquad \mbox{for all}\ t \in (0,\tm)   
  \eas
for some $c_1>0$ and, given any $p \in (1,\infty)$, we can pick $c_2(p)>0$ fulfilling
  \be{L^p-u-c_2}
   \io u^p(\cdot,t) \le c_2(p) \qquad \mbox{for all } t\in (0,\tm).  
  \ee
Hence, for any $q \in [1,\infty]$, by taking $p \in (1,\infty)$ such that 
$\frac{1}{p}-\frac{1}{q}<\frac{1}{2}$, we infer from the variation-of-constant formula and 
well-known smoothing estimates for the Neumann heat semigroup \cite[Lemma 1.3]{W2010_JDE} 
that  
  \bea{nabla-v-c_5}
  \|\nabla v(\cdot,t)\|_{L^q(\Omega)}
  &\le& c_3\|\nabla v_0\|_{L^q(\Omega)}
       +c_3\int_0^t \left( 1 + (t-s)^{-\frac{1}{2}-\left(\frac{1}{p}-\frac{1}{q}\right)}\right) e^{-\lambda (t-s)}
                      \|v(m_2-bu-v)\|_{L^p(\Omega)}\,ds \nn\\
  &\le& c_4\left(1+ c_1[(m_2+c_1)|\Omega|^\frac{1}{p}+bc_2(p)]\int_0^t  \left( 1 + (t-s)^{-\frac{1}{2}-\left(\frac{1}{p}-\frac{1}{q}\right)}\right) e^{-\lambda (t-s)}\,ds\right) \nn\\
  &\le& c_5  \qquad \mbox{for all } t\in (0,\tm)
  \eea
with some constants $c_3,c_4,c_5>0$ and $\lambda > 0$. 
In light of (\ref{L^p-u-c_2}) and (\ref{nabla-v-c_5}), we can deduce from  \cite[Lemma A.1]{TW2012_JDE} that 
  \bas
    \|u(\cdot,t)\|_{L^\infty} \le c_6 \qquad \mbox{for all}\ t \in (0,\tm)   
  \eas
for some $c_6>0$. Thus we have $\tm=\infty$ and the conclusion of Theorem \ref{theo_global} 
(i) holds. \qed
%
%
%
%
%
%
\subsection{Global existence of weak solutions in the case $n=3$}\label{sec: 3D}
In the case $n=3$, since we will not be able to obtain boundedness of 
$\|u(\cdot,t)\|_{L^\infty(\Omega)}$ in time, we shift to global classical 
solvablity of approximate problems (\ref{prob-general}) with 
$F(u):=\frac{u}{1+\eps u}$ and then pass to the limit $\eps \searrow 0$ 
as in the arguments in \cite{W2012} and \cite{W2016}.
More precisely, for $\eps \in (0,1)$ we consider 
\be{prob-ap}
    	\left\{ \begin{array}{rcll}
	& & \hspace*{-7mm}(\ueps)_t = \nabla\cdot ((d_1+\chi \veps)\nabla \ueps)
            - \chi \nabla \cdot \left(\dfrac{\ueps}{1+\eps \ueps}\nabla \veps\right) 
            +\ueps(m_1-\ueps+a\veps),  & x\in\Omega, \ t>0, \\[3mm]
	& & \hspace*{-7mm}(\veps)_t = d_2\Delta \veps +\veps(m_2-b\ueps-\veps),  & x\in\Omega, \ t>0, \\[2mm]
    & & \hspace*{-7mm}\frac{\pa \ueps}{\pa\nu} = \frac{\pa \veps}{\pa\nu} = 0,
	 & x\in\pO, \ t>0, \\[2mm]
	& & \hspace*{-7mm}
	\ueps(x,0)=u_0(x), \quad
	\veps(x,0)=v_0(x).
	 & x\in\Omega,   
	\end{array} \right.
\ee
As provided in Lemma \ref{lem_loc}, the approximate problem (\ref{prob-ap}) admits 
a unique local-in-time classical solution $(\ueps,\veps)$ up to a maximal time 
$\tmeps \in (0,\infty]$. 
Before proving that $\tmeps=\infty$  
we shall establish an energy inequality through a combination of $\io \ueps \log (\ueps)$ 
and $\io \frac{|\nabla \veps|^2}{\veps}$. 
\begin{lem}\label{lem_energy1}
Under the assumption of Lemma \ref{lem_loc}, there exists a constant $K>0$ such that 
the solution of (\ref{prob-ap}) satisfies
  \bea{energy1}
    & &\hspace{-6mm} \frac{d}{dt}\left[\io \ueps\log(\ueps)   
           + \frac{\chi}{2b} \io \frac{|\nabla \veps|^2}{\veps} \right] \nn\\ 
   & &\hspace{-6mm} 
        + \frac{1}{K}\left[
          \io \frac{d_1+\chi \veps}{\ueps}|\nabla \ueps|^2 
        + \io \ueps^2 \log (\ueps)
        + \io \frac{|D^2 \veps|^2}{\veps}
        + \io \frac{|\nabla \veps|^4}{\veps^3}\right] \nn\\
    & &\hspace{-6mm} \le K
    \qquad \mbox{for all}\ t \in (0,\tmeps).   
  \eea
\end{lem}
\proof%
From the first equation in (\ref{prob-ap}) and integration by parts it follows that 
  \bas
	& & \hspace{-13mm} 	\frac{d}{dt}\io \ueps \log (\ueps) \\
	&=& \io (\ueps)_t \log (\ueps) + \io (\ueps)_t \\
	&=& -\io \frac{d_1+\chi \veps}{\ueps}|\nabla \ueps|^2 
        + \chi \io \frac{1}{1+\eps\ueps}\nabla \ueps \cdot \nabla \veps  
        +\io \ueps(m_1-\ueps+a\veps)(\log (\ueps) +1). 
  \eas
Using (\ref{v-bound}) and noting that $\log (r)+1 \le r$ ($r>0$), we observe 
that if $\log (\ueps) +1 \ge 0$ then 
  \bas
   \ueps(m_1-\ueps+a\veps)(\log (\ueps) +1) \le \ueps(c_1-\ueps)(\log (\ueps) +1) 
   \le  - \frac{1}{2}\ueps^2\log (\ueps) + \frac{1}{2}c_1^2 \ueps 
  \eas
for some $c_1> 1$ and this holds even if $\log (\ueps) +1 < 0$ because 
  \bas
   \ueps(m_1-\ueps+a\veps)(\log (\ueps) +1) \le -\ueps \cdot \ueps \log (\ueps) 
   \le  \frac{1}{e} \ueps. 
  \eas
Hence by virtue of (\ref{mass}) we have
  \be{u-log}
	\frac{d}{dt}\io \ueps \log (\ueps)
	\le -\io \frac{d_1+\chi \veps}{\ueps}|\nabla \ueps|^2 
        + \chi \io \frac{1}{1+\eps\ueps} \nabla \ueps \cdot \nabla \veps  
        -\frac{1}{2}\io \ueps^2 \log (\ueps) +c_2  
  \ee
for some $c_2>0$. 
On the other hand, we infer from integration by parts and 
the second equation in (\ref{prob-ap}) that 
 \bas
	& & \hspace{-13mm} \frac{d}{dt}\io \frac{|\nabla \veps|^2}{\veps} \\ 
	&=& -2\io \frac{\Delta \veps}{\veps}(\veps)_t 
         + \io \frac{|\nabla \veps|^2}{\veps^2}(\veps)_t \\
	&=& -2d_2\io \frac{|\Delta \veps|^2}{\veps} 
        -2 \io \Delta \veps (m_2-b\ueps-\veps) 
        + d_2 \io \frac{|\nabla \veps|^2}{\veps^2}\Delta \veps 
        + \io |\nabla \veps|^2 \frac{m_2-b\ueps-\veps}{\veps} \\
	&\le& -2d_2\io \frac{|\Delta \veps|^2}{\veps} 
        -2b \io \nabla \ueps \cdot \nabla \veps  
        + d_2 \io \frac{|\nabla \veps|^2}{\veps^2}\Delta \veps  
        + m_2 \io \frac{|\nabla \veps|^2}{\veps} - 3\io |\nabla \veps|^2
  \eas
for all $t \in (0,\tmeps)$.
Here, thanks to \cite[Lemma 3.1]{W2012}, we have
  \bas
      -\io \frac{|\Delta \veps|^2}{\veps} 
    = -\io \frac{|D^2 \veps|^2}{\veps} 
      -\frac{3}{2} \io \frac{|\nabla \veps|^2}{\veps^2}\Delta \veps 
      +\io \frac{|\nabla \veps|^4}{\veps^3}
      +\frac{1}{2} \int_{\pO} \frac{1}{\veps} \frac{\pa}{\pa \nu}|\nabla \veps|^2
  \eas
and as in \cite[Proof of Lemma 3.2]{W2012} we can see from 
a direct computation and integration by parts that
  \bas
      \io \veps |D^2\log (\veps)|^2 
    = \io \frac{|D^2 \veps|^2}{\veps} + \io \frac{|\nabla \veps|^2}{\veps^2}\Delta \veps 
     -\io \frac{|\nabla \veps|^4}{\veps^3},
  \eas
whence adding these two identities entail
  \bas
      -\io \frac{|\Delta \veps|^2}{\veps}
      +\io \veps |D^2\log(\veps)|^2
    = -\frac{1}{2} \io \frac{|\nabla \veps|^2}{\veps^2}\Delta \veps 
      +\frac{1}{2} \int_{\pO} \frac{1}{\veps} \frac{\pa}{\pa \nu}|\nabla \veps|^2.
  \eas
Thus we obtain
 \be{v-energy}
	\frac{d}{dt}\io \frac{|\nabla \veps|^2}{\veps}  
	\le -2d_2 \io \veps |D^2\log(\veps)|^2 
        + d_2 \int_{\pO} \frac{1}{\veps} \frac{\pa}{\pa \nu}|\nabla \veps|^2 
        -2b \io \nabla \ueps \cdot \nabla \veps 
        + m_2 \io \frac{|\nabla \veps|^2}{\veps}
  \ee
for all $t \in (0,\tmeps)$. 
Therefore, adding (\ref{u-log}) and (\ref{v-energy}) multiplied by $\frac{\chi}{2b}$, 
we have that with some constants $c_3,c_4>0$, 
  \bea{energy00}
    & &\hspace{-6mm} \frac{d}{dt}\left[\io \ueps\log (\ueps) + \frac{\chi}{2b} \io \frac{|\nabla \veps|^2}{\veps} \right] \nn\\
    & & + \io \frac{d_1+\chi \veps}{\ueps}|\nabla \ueps|^2 
        + \frac{1}{2}\io \ueps^2 \log (\ueps)
        + \frac{\chi d_2}{b} \io \veps |D^2\log (\veps)|^2 \nn\\
    & &\hspace{-9mm} \le c_2+ c_3 \int_{\pO} \frac{1}{\veps} \frac{\pa}{\pa \nu}|\nabla \veps|^2 
           +c_4 \io \frac{|\nabla \veps|^2}{\veps}
           -\chi\eps \io \frac{\ueps}{1+\eps\ueps}\nabla \ueps \cdot \nabla \veps
  \eea
for all $t \in (0,\tmeps)$.  
Invoking \cite[(3.9)]{W2016}, we can estimate the last term on the left-hand side 
of (\ref{energy00}) as  
  \bas
       \frac{\chi d_2}{b} \io \veps |D^2\log (\veps)|^2 
   \ge c_5 \io \frac{|D^2 \veps|^2}{\veps}
     + c_5 \io \frac{|\nabla \veps|^4}{\veps^3}
  \eas
for some $c_5>0$, while Young's inequality and (\ref{v-bound}) give  
  \be{Young-v}
    c_4 \io \frac{|\nabla \veps|^2}{\veps} 
   \le \frac{c_5}{8}\io \frac{|\nabla \veps|^4}{\veps^3} +c_6
  \ee
as well as 
  \bas
    -\chi\eps \io \frac{\ueps}{1+\eps\ueps}\nabla \ueps \cdot \nabla \veps
    = \io \frac{\nabla \ueps }{\ueps^\frac{1}{2}}\cdot
      \frac{\nabla \veps}{\veps^\frac{3}{4}}\cdot
      \frac{-\chi\eps \ueps^\frac{3}{2}\veps^\frac{3}{4}}{1+\eps\ueps} 
  \le \frac{d_1}{2}\io \frac{|\nabla \ueps|^2}{\ueps}
       +\frac{c_5}{8}\io \frac{|\nabla \veps|^4}{\veps^3} +c_7\io \ueps^2
  \eas
for some $\eps$-independent constants $c_6, c_7>0$, 
all of which imply that 
  \bea{energy01}
    & &\hspace{-6mm} \frac{d}{dt}\left[\io \ueps\log (\ueps) + \frac{\chi}{2b} \io \frac{|\nabla \veps|^2}{\veps} \right] \nn\\
    & & + \frac{1}{2}\io \frac{d_1+\chi \veps}{\ueps}|\nabla \ueps|^2 
        + \frac{1}{2}\io \ueps^2 \log(\ueps)
        +c_5 \io \frac{|D^2 \veps|^2}{\veps}
        +\frac{3}{4}c_5 \io \frac{|\nabla \veps|^4}{\veps^3} \nn\\
    & &\hspace{-9mm} \le c_2 + c_3 \int_{\pO} \frac{1}{\veps} \frac{\pa}{\pa \nu}|\nabla \veps|^2 
           +c_6+c_7\io \ueps^2
  \eea
for all $t \in (0,\tmeps)$. 
In the same way as in the proof of \cite[(3.10)]{ISY}, 
noting that $\frac{\pa}{\pa \nu}|\nabla \veps|^2  \le c_8 |\nabla \veps|^2$ 
for some $c_8>0$ and combining the Sobolev compact embedding 
$W^{r+\frac{1}{2},2}(\Omega) 
  \hookrightarrow L^2(\pa\Omega)$ ($r \in (0,\frac{1}{2})$)
with the fractional Gagliardo-Nirenberg inequality, we can obtain
  \bas
      \int_{\pO} \frac{1}{\veps} \frac{\pa}{\pa \nu}|\nabla \veps|^2 
    \le 
        4c_8 \|\nabla \sqrt{\veps}\|_{L^2(\pO)}^2 
    \le c_9 \|D^2\sqrt{\veps}\|_{L^2(\Omega)}^{2\alpha}
            \|\nabla \sqrt{\veps}\|_{L^2(\Omega)}^{2(1-\alpha)}
      +c_9\|\nabla \sqrt{\veps}\|_{L^2(\Omega)}^{2}
  \eas
for some $c_9>0$ and $\alpha \in (0,1)$. Here a direct computation gives
  \bas
     \|D^2\sqrt{\veps}\|_{L^2(\Omega)}^2
     \le c_{10}\left( \io \frac{|D^2 \veps|^2}{\veps}+ \io \frac{|\nabla \veps|^4}{\veps^3}\right)
  \eas
with some $c_{10}>0$, 
and like in (\ref{Young-v}) for any small $\delta>0$ 
we can take $c_{11}(\delta)>0$ such that 
  \bas
     \|\nabla \sqrt{\veps}\|_{L^2(\Omega)}^{2}
     = \frac{1}{4}\io \frac{|\nabla \veps|^2}{\veps} 
     \le \delta \io \frac{|\nabla \veps|^4}{\veps^3} +c_{11}(\delta).
  \eas
Therefore the Young inequality yields
  \bas
       c_3\int_{\pO} \frac{1}{\veps} \frac{\pa}{\pa \nu}|\nabla \veps|^2 
   \le 
       \frac{c_5}{2} \io \frac{|D^2 \veps|^2}{\veps}
        +\frac{c_5}{4} \io \frac{|\nabla \veps|^4}{\veps^3}+c_{12}
  \eas
for some $c_{12}>0$. 
Plugging this inequality into (\ref{energy01}) and employing the inequality 
$c_7 r^2 \le \frac{1}{4}r^2 \log (r) +c_{13}$ $(r \ge 0)$ with some $c_{13}>0$ in 
the last summand in  (\ref{energy01}) imply 
the desired inequality (\ref{energy1}). 
\qed
As a consequence of Lemma \ref{lem_energy1} 
we can draw the following a priori estimates. 
\begin{lem}\label{lem_energy2}
  Under the assumption of Lemma \ref{lem_loc}, assume further that 
 $u_0 \in L\log L(\Omega)$ and $\sqrt{v_0} \in W^{1,2}(\Omega)$. 
  Then there exists $C>0$ such that for all $\eps \in (0,1)$  
  the solution of (\ref{prob-ap}) satisfies
  \be{energy-pointwise}
    \io \ueps(\cdot,t)\log (\ueps(\cdot,t))
    + \frac{\chi}{2b} \io \frac{|\nabla \veps(\cdot,t)|^2}{\veps(\cdot,t)} 
    \le C 
    \qquad \mbox{for all}\ t \in (0,\tmeps)   
  \ee
and 
  \be{energy-integral}
     \int_0^T \io \frac{d_1+\chi \veps}{\ueps}|\nabla \ueps|^2 
   + \int_0^T \io \ueps^2 \log (\ueps)
   + \int_0^T \io \frac{|D^2 \veps|^2}{\veps}
   + \int_0^T \io \frac{|\nabla \veps|^4}{\veps^3}
   \le C\cdot(T+1)
  \ee
for all $T \in (0,\tmeps)$. 
\end{lem}
\proof 
We first make use of the trivial inequality 
$r\log (r) \le r^2\log (r) $ $(r \ge 0)$ to 
obtain 
  \bas
       \io \ueps \log (\ueps) \le \io \ueps^2 \log  (\ueps). 
  \eas
On the other hand, we next invoke the Young inequality and (\ref{v-bound}) to infer
  \bas
       \frac{\chi}{2b} \io \frac{|\nabla \veps|^2}{\veps} 
   \le \io \frac{|\nabla \veps|^4}{\veps^3} +c_1
  \eas
for some $c_1>0$. Using these two inequality in (\ref{energy1}), we observe that 
  \bas
    y_\eps(t):=\io \ueps(\cdot,t)\log (\ueps(\cdot,t)) 
    + \frac{\chi}{2b} \io \frac{|\nabla \veps(\cdot,t)|^2}{\veps(\cdot,t)} 
  \eas
satisfies
  \be{ODI-uniform}
    y_\eps'(t) + \frac{1}{K}y_\eps(t) 
    \le c_2:=K + \frac{c_1}{K}
    \qquad \mbox{for all}\ t \in (0,\tmeps).   
  \ee
This provides 
   \bas
     y_\eps(t) \le c_3:=\max\left\{\io u_0 \log (u_0)
     + \frac{\chi}{2b} \io \frac{|\nabla v_0|^2}{v_0},\ 
       Kc_2\right\} 
   \eas  
for all $t \in (0,\tmeps)$, which precisely warrants (\ref{energy-pointwise}). 
Integrating (\ref{energy1}), we arrive at (\ref{energy-integral}).
\qed
\begin{lem}\label{lem_global-app}
  Under the assumption of Lemma \ref{lem_energy2}, 
  for all $\eps \in (0,1)$, we have $\tmeps=\infty$; that is, 
  the solution of (\ref{prob-ap}) is global in time. 
\end{lem}
\proof 
To facilitate a proof by contradiction, we assume that $\tmeps < \infty$.
Then, multiplying the first equation in (\ref{prob-ap}) by $\ueps^3$ and 
noting that $\frac{\ueps}{1+\eps\ueps} \le \frac{1}{\eps}$,  
we see from integration by parts, (\ref{v-bound}) and Young's inequality that
  \bas
       \frac{1}{4}\frac{d}{dt}\io \ueps^4 
   &=& -3 \io (d_1+\chi \veps)\ueps^2|\nabla \ueps|^2
       +3\chi \io \frac{\ueps^3}{1+\eps\ueps}\nabla \ueps \cdot \nabla \veps 
       +\io \ueps^4(m_1-\ueps+a\veps) \\
   &\le& -3d_1 \io \ueps^2|\nabla \ueps|^2
       +c_1 \io \ueps |\nabla \ueps| \cdot \ueps \cdot 
         \frac{|\nabla \veps|}{\veps^\frac{3}{4}}
       +c_2\io \ueps^4 \\
   &\le& c_3 \io \frac{|\nabla \veps|^4}{\veps^3} 
       +c_4 \io \ueps^4 
  \eas
for some $c_1, c_2, c_3, c_4>0$. In view of (\ref{energy-integral}) this enables us to 
find $c_5>0$ fulfilling
  \be{L^4-u}
    \io \ueps^4(\cdot,t) \le c_5  \qquad \mbox{for all}\ t \in (0,\tmeps).
  \ee
Hence, similarly in (\ref{nabla-v-c_5}),  we have 
  \bea{nabla-v-bdd}
  \|\nabla v(\cdot,t)\|_{L^\infty(\Omega)}
  &\le& c_6\|\nabla v_0\|_{L^\infty(\Omega)}
       +c_6\int_0^t \left( 1 + (t-s)^{-\frac{1}{2}-\frac{3}{2}\cdot\frac{1}{4}} \right)e^{-\lambda (t-s)}
                     \|v(m_2-bu-v)\|_{L^4(\Omega)}\,ds \nn\\
  &\le& c_7 
  \qquad \mbox{for all } t\in (0,\tmeps)
  \eea
with some constants $c_6,c_7>0$ and $\lambda > 0$. 
Therefore, for all $p \in (1,\infty)$, 
again by multiplying the first equation in (\ref{prob-ap}) by $\ueps^{p-1}$, 
we obtain
  \bas
       \frac{1}{p}\frac{d}{dt}\io \ueps^p 
   &\le & -(p-1)d_1 \io \ueps^{p-2}|\nabla \ueps|^2
       +(p-1)\chi \io \frac{\ueps^{p-1}}{1+\eps\ueps}\nabla \ueps \cdot \nabla \veps 
       +\io \ueps^{p}(m_1-\ueps+a\veps) \\
   &\le& -(p-1)d_1 \io \ueps^{p-2}|\nabla \ueps|^2
       +c_8 \io \ueps^{p-1} |\nabla \ueps|  
       +c_9\io \ueps^p \\
   &\le&  c_{10}\io \ueps^p
  \eas
for some $c_8, c_9, c_{10}>0$. This implies that for all $p \in (1,\infty)$ 
there is $c_{11}(p)>0$ such that
  \bas
    \io \ueps^p(\cdot,t) \le c_{11}(p)  \qquad \mbox{for all } t\in (0,\tmeps), 
  \eas
so that an iteration argument as in \cite[Lemma A.1]{TW2012_JDE} yields 
the existence of $c_{12}>0$ satisfying
  \bas
     \|\ueps(\cdot,t)\|_{L^\infty(\Omega)} \le c_{12} 
     \qquad \mbox{for all } t\in (0,\tmeps). 
  \eas
Combined with (\ref{nabla-v-bdd}), this means that $\tmeps=\infty$ as our initial assumption combined with the aforementioned bounds leads to a contradiction according to the extensibility criterion (\ref{ext}).  
\qed
The following further estimates for spatio-temporal integrals are 
immediate consequences of Lemmas \ref{lem_basic-u} and \ref{lem_energy2}.   
\begin{cor}\label{lem_further-estimates}
  Under the assumption of Lemma \ref{lem_energy2}, 
  there exists $C>0$ such that for all $\eps \in (0,1)$  
  the solution of (\ref{prob-ap}) satisfies 
  \be{int-u2-estimate}
   \int_0^T \io u_\eps^{2} \le C\cdot(T+1)   \qquad \mbox{for all}\ T>0
  \ee
and 
  \be{int-nabla-u-estimate}
   \int_0^T \io  |\nabla u_\eps|^{\frac{4}{3}} \le C\cdot(T+1)   \qquad \mbox{for all}\ T>0
  \ee
as well as
  \be{int-nabla-v-estimate}
   \int_0^T \io |\nabla v_\eps|^4 \le C\cdot(T+1)   \qquad \mbox{for all}\ T>0.
  \ee
\end{cor}
\proof 
In view of Lemma \ref{lem_basic-u} we find $c_1>0$ fulfilling  
  \bas
	\int_0^T \io \ueps^2 \le c_1\cdot(T+1) \qquad \mbox{for all } T>0, 
  \eas
which means (\ref{int-u2-estimate}). According to Lemma \ref{lem_energy2}, 
there exists $c_2>0$ such that
  \be{int-estimates-com}
     \int_0^T \io \frac{|\nabla \ueps|^2}{\ueps}
    +\int_0^T \io \frac{|\nabla \veps|^4}{\veps^3}
   \le c_2 \cdot(T+1) \qquad \mbox{for all } T>0.
  \ee
Thus, the H\"older inequality gives
  \bas
        \int_0^T \io |\nabla u_\eps|^{\frac{4}{3}}
    &=& \int_0^T \io \ueps^\frac{2}{3}\cdot
          \frac{|\nabla u_\eps|^{\frac{4}{3}}}{\ueps^\frac{2}{3}} \\
    &\le& \left(\int_0^T \io \ueps^2 \right)^{\frac{1}{3}}
          \left(\int_0^T \io \frac{|\nabla \ueps|^2}{\ueps}\right)^{\frac{2}{3}} \\
    &\le& c_1^{\frac{1}{3}}c_2^{\frac{2}{3}}(T+1) \qquad \mbox{for all}\ T>0, 
  \eas
which shows (\ref{int-nabla-u-estimate}). 
Finally, (\ref{int-estimates-com}) together with (\ref{v-bound}) 
implies (\ref{int-nabla-v-estimate}).  
\qed

In preparation for the shortly following compactness argument used as the source for our weak solutions, we now derive some additional uniform estimates for the time derivatives of our approximate solutions.

\begin{lem}\label{lem_time-derivative-estimates}
  Under the assumption of Lemma \ref{lem_energy2}, 
  there exists $C>0$ such that for all $\eps \in (0,1)$  
  the solution of (\ref{prob-ap}) satisfies 
\be{int-ut-estimate}
	\int_0^T \|u_{\eps t}(\cdot, t)\|_{(W^{1,4}(\Omega))^*} \, d t \le C\cdot(T+1)   \qquad \mbox{for all}\ T>0
\ee
and
\be{int-vt-estimate}
	\int_0^T \|v_{\eps t}(\cdot, t)\|_{(W^{1,4}(\Omega))^*} \, d t \le C\cdot(T+1)   \qquad \mbox{for all}\ T>0.
\ee
\end{lem} 
\proof%
Let $\varphi$ be a function in $W^{1,4}(\Omega) \subset L^\infty(\Omega)$. We then test the first equation in (\ref{prob-ap}) with $\varphi$ and apply 
integration by parts as well as the Young and H\"older inequalities to see that there exists $c_1 > 0$ such that
\bas
	 \left|\int_\Omega u_{\eps t} \varphi \right| 
	&\leq& \int_\Omega  (d_1 + \chi \veps )|\grad \ueps||\grad \varphi|+ \chi \int_\Omega \frac{\ueps}{1+\eps\ueps}|\grad \veps||\grad \varphi| + \int_\Omega |\varphi| \ueps(m_1 + \ueps + a \veps) \\
	&\leq& (d_1 + \chi \|v_\eps\|_{L^\infty(\Omega)})\|\grad \ueps\|_{L^\frac{4}{3}(\Omega)}\| \grad \varphi \|_{L^4(\Omega)} + \chi \|\ueps\|_{L^2(\Omega)}\|\grad \veps\|_{L^4(\Omega)}\|\grad \varphi\|_{L^4(\Omega)} \\
	& & + \|\varphi\|_{L^\infty(\Omega)} \int_\Omega \ueps(m_1 + \ueps + a \veps)    \\
	&\leq& c_1 \left( \int_\Omega |\grad \ueps|^\frac{4}{3} + \int_\Omega \ueps^2 +  \int_\Omega |\grad \veps|^4 + \int_\Omega \veps^2 +  \|\veps\|^4_{L^\infty(\Omega)} + 1\right) \left( \| \grad \varphi \|_{L^4(\Omega)}  + \|\varphi\|_{L^\infty(\Omega)} \right)
\eas
for all $t \in (0,\infty)$. Due to the bounds established in Corollary \ref{lem_further-estimates} and Lemma \ref{lem_basic-v} as well as the fact that in three dimensions $W^{1,4}(\Omega)$ is continuously embedded into $L^\infty(\Omega)$, we immediately gain the bound (\ref{int-ut-estimate}) from the above. By a similar testing procedure applied to the second equation in (\ref{prob-ap}), we gain $c_2 > 0$ such that
\bas
	\left|\int_\Omega v_{\eps t} \varphi \right| &\leq& d_2 \int_\Omega |\grad \veps| |\grad \varphi | + \int_\Omega |\varphi| \veps (|m_2| + b \ueps + \veps) \\
	&\leq & c_2 \left( \int_\Omega \ueps^2 + \int_\Omega |\grad \veps|^4  + \int_\Omega \veps^2  + 1 \right)\left( \|\grad \varphi \|_{L^\frac{4}{3}(\Omega)} + \|\varphi\|_{L^\infty(\Omega)} \right)
\eas
for all $t \in (0,\infty)$. Again by application of Corollary \ref{lem_further-estimates} and Lemma \ref{lem_basic-v}, this implies the second bound (\ref{int-vt-estimate}) and therefore completes the proof.
\qed
We now construct candidates for our weak solutions as limits of our approximate solutions along a suitable null sequence $(\eps_j)_{j\in\N}$ by applying the Aubin--Lions compact embedding theorem as well as the weak compactness properties of bounded sets in Sobolev spaces.
\begin{lem}\label{lem_limit-process}
	Under the assumption of Lemma \ref{lem_energy2}, there exist a null sequence $(\eps_j)_{j\in\N}$ and a.e.\ nonnegative functions $u, v: \overline{\Omega} \times [0,\infty) \rightarrow \R$ such that $v \in L^\infty(\Omega \times (0,\infty))$ as well as
	\bea{}
		\ueps &\rightarrow& u \qquad \text{ in } \qquad L^2_{\rm loc}([0,\infty); L^2(\Omega)) \text{ and  a.e.\ in } \Omega\times(0,\infty), \label{ueps_strong_convergence} \\
		\ueps &\rightharpoonup& u \qquad \text{ in } \qquad L^\frac{4}{3}_{\rm loc}([0,\infty); W^{1,\frac{4}{3}}(\Omega)), \label{ueps_weak_convergence} \\
		\veps &\rightarrow& v \qquad \text{ in } \qquad L^4_{\rm loc}([0,\infty); L^4(\Omega)) \text{ and  a.e.\ in } \Omega\times(0,\infty) \text{ and } \label{veps_strong_convergence}  \\
		\veps &\rightharpoonup& v \qquad \text{ in } \qquad L^4_{\rm loc}([0,\infty); W^{1,4}(\Omega)) \label{veps_weak_convergence} 
	\eea
	as $\eps = \eps_j \searrow 0$.
\end{lem}
\proof%
According to Corollary \ref{lem_further-estimates}, the family $(\ueps)_{\eps \in (0,1)}$ is bounded in  $L^\frac{4}{3}((0,T); W^{1,\frac{4}{3}}(\Omega))$ for all $T > 0$ and, according to Lemma \ref{lem_time-derivative-estimates}, the family $(u_{\eps t})_{\eps \in (0,1)}$ is bounded in  $L^1((0,T); (W^{1,4}(\Omega))^*)$ for all $T > 0$. Thus applying the Aubin--Lions lemma to the triple of spaces $W^{1,\frac{4}{3}}(\Omega) \subset\subset L^\frac{4}{3}(\Omega)\subset(W^{1,4}(\Omega))^*$ as well as the weak compactness property of bounded sets in Sobolev spaces for all $T\in\N$ combined with a diagonal sequence argument yields a null sequence $(\eps_j)_{j\in\N}$ and $u: \overline{\Omega} \times [0,\infty) \rightarrow \R$ such that (\ref{ueps_weak_convergence}) holds and such that $\ueps \rightarrow u$ in $L^{\frac{4}{3}}_{\rm loc}([0,\infty); L^\frac{4}{3}(\Omega))$ as well as a.e.\ on $\Omega\times(0,\infty)$ as $\eps = \eps_j \searrow 0$. As further Lemma \ref{lem_energy1} provides $c > 0$ such that
\bas
     \int_0^T \int_\Omega \ueps^2 \log(\ueps) \leq c \cdot (T + 1)    \qquad \mbox{for all}\ T>0,
\eas
a combination of the De La Vallée Poussin criterion for uniform integratiblity and the Vitali convergence theorem yields the convergence property (\ref{ueps_strong_convergence}).
\abs
Similarly according to Lemma \ref{lem_basic-v} and Corollary \ref{lem_further-estimates}, the family $(\veps)_{\eps \in (0,1)}$ is bounded in the space $L^4((0,T); W^{1,4}(\Omega))$ for all $T > 0$ and, according to Lemma \ref{lem_time-derivative-estimates}, the family $(v_{\eps t})_{\eps \in (0,1)}$ is bounded in the space $L^1((0,T); (W^{1,4}(\Omega))^*)$ for all $T > 0$. Given this, we can again apply the Aubin--Lions lemma to the triple of spaces $W^{1,4}(\Omega) \subset \subset L^4(\Omega) \subset (W^{1,4}(\Omega))^*$ and use the same weak compactness properties of bounded sets in Sobolev spaces to immediately gain a function $v:\overline{\Omega} \times [0,\infty) \rightarrow \R$ such that both (\ref{veps_strong_convergence}) and (\ref{veps_weak_convergence}) hold (after potentially switching to a subsequence).
\abs
The a.e.\ nonnegativity of both $u$ and $v$ as well as the global boundedness of $v$ then follows from the positivity of the approximate solutions and the boundedness properties laid out in Lemma \ref{lem_basic-v} combined with the a.e.\ convergence already established in (\ref{ueps_strong_convergence}) and (\ref{veps_strong_convergence}).
\qed
As we have now constructed our solution candidates, it remains only to be shown that they in fact fulfill our desired weak solution property and thus prove the second half of Theorem \ref{theo_global}. \abs
\proofc of Theorem \ref{theo_global} (ii).\quad Let $u$ and $v$ as well as the null sequence $(\eps_j)_{j\in\N}$ be as constructed in Lemma \ref{lem_limit-process}. Notably given the convergence, nonnegativity and boundedness properties established in said lemma, we immediately gain the necessary regularity properties for our weak solution definition. We thus only need to further show that $u$ and $v$ fulfill (\ref{weak_solution_u}) and (\ref{weak_solution_v}). To this end, we first note that each $\ueps$ already fulfills a slightly modified version of (\ref{weak_solution_u}) corresponding to the first equation in (\ref{prob-ap}) and each $\veps$ already fulfills (\ref{weak_solution_v}) exactly. As such, it now only remains to show that all terms in the weak solution equation solved by our approximate solutions converge to their counterparts without $\eps$ as $\eps = \eps_j \searrow 0$ as well as handle the slightly different taxis term. The former can be easily seen by application of the convergence properties from Lemma \ref{lem_limit-process} while we will now discuss the latter in a little more detail. To facilitate this, we first fix $\varphi \in C^\infty_0(\overline{\Omega}\times[0,\infty))$. Due to the dominated convergence theorem and the a.e.\ pointwise convergence given in (\ref{ueps_strong_convergence}) we can conclude that $\frac{\grad \varphi}{1 + \eps\ueps}$ converges to $\grad \varphi$ in $L^4(\Omega\times(0,\infty))$ as $\eps = \eps_j \searrow 0$. Combining this with (\ref{ueps_strong_convergence}) then yields that $\frac{\ueps \grad \varphi}{1 + \eps\ueps}$ converges to $u \grad \varphi$ in $L^\frac{4}{3}(\Omega\times(0,\infty))$, which in turn combined with (\ref{veps_weak_convergence}) implies that $\int_0^\infty\int_\Omega \ueps \grad \veps \cdot \frac{\grad \varphi}{1 + \eps \ueps}$ converges to $\int_0^\infty\int_\Omega u \grad v \cdot \grad \varphi$ as $\eps = \eps_j \searrow 0$ and therefore completes the proof.
\qed
%
%
%
%
%
%
%
%
\mysection{Stabilization}\label{sec: stabilize}

As we have now proven the first of our two main results, we shift our focus from the aforementioned existence result to our analysis of long-time behavior. Similarly to the previous section, we begin by deriving some a priori bounds for the solutions to (\ref{prob-general}). Importantly, we again take great care to ensure that all the constants involved in this process are independent of $F$ to enable us to not only use said bounds in the two-dimensional case but also in the three-dimensional case, where the solution construction relies on approximation.
\abs
In an effort to eliminate an otherwise necessary initial data condition on $v_0$, we first show that after an appropriate waiting time the second component of solutions to (\ref{prob-general}) approaches $m_{2+} = \max(0, m_2)$ uniformly from above. 
\begin{lem}\label{lem_eventual_v_bound}
	Under the assumption of Lemma \ref{lem_loc} and for each $m > m_{2+} = \max(0, m_2)$, there exists $T \equiv T(m) > 0$ such that
	\[
	v(x, t) \leq m
	\]
	for all $x\in\overline{\Omega}$ and $t \in (T, \infty)$.
\end{lem}
\proof%
By comparison with the solution to the initial value problem $y' = y(m_2-y)$, $y(0) = \|v_0\|_{L^\infty(\Omega)}$, we gain that
\bas
	v(\cdot, t) \leq \begin{cases} 		
		\left(\frac{1}{m_2} + \left(\frac{1}{\|v_0\|_{L^\infty(\Omega)}} - \frac{1}{m_2}\right)e^{-m_2t}\right)^{-1} & \text{ if } m_2 \neq 0 \\
		\left( \frac{1}{\|v_0\|_{L^\infty(\Omega)}} +  t \right)^{-1} &\text{ if } m_2 = 0
	\end{cases}
\eas
for all $t \in (0,\infty)$. Our desired result then follows directly from this as the right-hand side converges to $m_2$ if $m_2 > 0$ and to $0$ if $m_2 \leq 0$ as $t \rightarrow \infty$.
\qed
%
%
%
%
%
%
\subsection{An eventual uniform energy-type inequality for solutions to (\ref{prob-general})}
As is not uncommon (see e.g.\ \cite{W2020_ANS}), our argument for the stabilization behavior laid out in Theorem \ref{theo_stabilization} chiefly rests on the derivation of an eventual (uniform) energy-type inequality for solutions of (\ref{prob-general}), which in the two-dimensional case can then be directly applied to the classical solutions we have constructed in the previous section or whose most important consequences survive the limit process used in our construction of weak solutions in the three-dimensional case. Notably, the energy functional involved is different and in a sense weaker from the one employed in Section \ref{sec: 3D} as our long-time stabilization argument in turn necessitates a stronger type of energy inequality than the aforementioned existence result.
\abs
To establish an important building block for the functionals involved in said energy-type inequality, we now define the family of nonnegative functions
\[
H_\xi(\eta) := \begin{cases} 
\eta - \xi - \xi \log\left(\frac{\eta}{\xi}\right) &\text{ if } \xi > 0,\\
\eta &\text{ if } \xi = 0  
\end{cases}
\]
for all $\xi \geq 0, \eta \geq 0$
and recall the following functional inequalities involving the aforementioned functions as e.g.\ established in \cite{W2020_ANS}.
\begin{lem}\label{lem_H_props}
	Let $\Omega$ be a smooth bounded domain and  $\xi \geq 0$. Then
	\be{Hxi_lower_bound}
	\int_\Omega |\varphi - \xi| \leq \frac{1}{1-\log(2)} \int_\Omega H_\xi(\varphi) + \sqrt{8\xi |\Omega|} \left\{ \int_\Omega H_\xi(\varphi) \right\}^\frac{1}{2}
	\ee for all measurable $\varphi : \Omega \rightarrow (0,\infty)$ and there exists $C > 0$ such that
	\be{Hxi_upper_bound}
	\int_\Omega H_\xi^2(\varphi) \leq \int_\Omega \varphi^2 + C \xi^2 \left( \int_\Omega \frac{|\grad \varphi|^2}{\varphi^2} +\left\{\int_\Omega |\log \varphi|\right\}^2 + 1 \right)
	\ee
	for all positive $\varphi \in L^2(\Omega)$ with $\log \varphi \in W^{1,2}(\Omega)$.
\end{lem}
\proof%
	The case $\xi = 0$ is trivial while the remaining case of $\xi > 0$ follows from \cite[Lemma 3.1, Lemma 3.2]{W2020_ANS}.
\qed
Using the now defined functions $H_{\xi}(\eta)$ and assuming that (\ref{stablization_param_condition}) holds, we can define the energy-type functional
\[
	E(u,v) := \int_\Omega H_\us(u) + \frac{a}{b}\int_\Omega H_\vs(v) + \frac{2}{b^2 \widehat{m_2}}\int_\Omega (v - \vs)^2
\]
for all measurable functions $u,v: \Omega \rightarrow [0,\infty]$ as well as the functional
\[
	G(u,v) := \int_\Omega \frac{|\grad u|^2}{u^2} + \int_\Omega \frac{|\grad v|^2}{v^2} + \int_\Omega (u - \us)^2 + \int_\Omega (v - \vs)^2 
\]
for all $u,v \in L^2(\Omega)$ with $\log(u), \log(v) \in W^{1,2}(\Omega)$.
Here $\us$ and $\vs$ are defined as in (\ref{equilibria}) and $\widehat{m_2} > m_{2+} = \max(0, m_2)$ is chosen in such a way that
\be{approx_stablization_param_condition}
	\chi^2 < \frac{4 d_1 d_2}{b\widehat{m_2}\us} \left( \frac{a \vs}{\widehat{m_2}} +\frac{4}{b} \right).
\ee  
Having introduced the central functionals of this section, we can now derive our desired eventual energy inequality as well as ensure that the energy $E$ remains uniformly bounded in the time interval before the energy inequality provides us with even stronger guarantees. The latter is especially important for the three-dimensional case. Note that, here as well as later in this section, we will make extensive use of the fact that $\us$ and $\vs$ are chosen to exactly fulfill $(m_1 - \us + a\vs) = 0$ and $(m_2 - b\us - \vs) = 0$ if $\vs > 0$ or $(m_2 - b\us - \vs) \leq 0$ if $\vs = 0$. 
\begin{lem}\label{lem_stab_energy}
	Under the assumption of Lemma \ref{lem_loc} and if (\ref{stabilization_initial_data_condition}) as well as (\ref{stablization_param_condition}) hold, there exist $C > 0$, $T_E > 0$ and $\delta > 0$ such that
	\be{stab_energy_bound}
		E(u(\cdot, t), v(\cdot, t)) \leq C
	\ee
	for all $t \in [0,T_E]$ and
	\be{stab_energy_ineq}
	\frac{d}{dt} E(u(\cdot, t), v(\cdot, t)) \leq - \delta\, G(u(\cdot, t), v(\cdot, t))
	\ee
	for all $t \in (T_E,\infty)$.
\end{lem}

\proof%
We begin by fixing $\delta > 0$ such that
\be{energy_param_1}
	\delta \leq \frac{a}{b}
\ee
and
\be{energy_param_2}
	(1 - \delta) \widehat{m_2} > m_{2+} = \max(0, m_2).
\ee
as well as
\bas
	\chi^2 < \frac{4 (d_1 - \frac{\delta}{\us}) d_2}{b \widehat{m_2} \us } \left( \frac{a \vs}{\widehat{m_2}} +\frac{4}{b} - \delta\frac{b}{d_2 \widehat{m_2}} \right)
\eas  
while ensuring that all factors on the right-hand side stay positive individually, which is possible due to (\ref{approx_stablization_param_condition}). This then implies 
\be{energy_param_3}
	\left(\frac{\chi^2 \us}{4 (d_1 - \frac{\delta}{\us})} - \frac{4d_2}{b^2 \widehat{m_2}}  \right) \widehat{m_2}^2 - d_2 \vs \frac{a}{b} < - \delta
\ee 
after some slight rearrangement.
As $\delta$ fulfills (\ref{energy_param_2}), we can use Lemma \ref{lem_eventual_v_bound} to fix $T_E > 0$ such that 
\bea{energy_eventual_m_bound}
	\|v(\cdot, t)\|_{L^\infty(\Omega)} \leq (1-\delta) \widehat{m_2} \leq \widehat{m_2}
\eea
for all $t \in (T_E, \infty)$.
\abs
Now integrating the first equation in (\ref{prob-general}) over $\Omega$ yields 
\bas
	\frac{d}{dt}\int_\Omega u &=& \int_\Omega u(m_1 - u + a v)
\eas
for all $t\in(0,\infty)$ due to the Neumann boundary conditions.
Testing the same equation with $-\frac{1}{u}$ further gives us
\bas
	-\frac{d}{d t}\int_\Omega \log(u) &=& -\int_\Omega \frac{u_t}{u} = -\int_\Omega \frac{\div ((d_1+\chi v) \grad u)}{u} + \chi \int_\Omega \frac{1}{u} \div \left( F(u) \grad v \right) - \int_\Omega (m_1 - u + a v)  \\
	&=& -\int_\Omega (d_1 + \chi v)\frac{|\grad u|^2}{u^2} + \chi \int_\Omega \frac{F(u)}{u^2} \grad u \cdot \grad v - \int_\Omega (m_1 - u + a v) \\
	&\leq& -\frac{\delta}{\us} \int_\Omega \frac{|\grad u|^2}{u^2} + \frac{\chi^2}{4  (d_1 - \frac{\delta}{\us})} \int_\Omega \left(\frac{F(u)}{u}\right)^2|\grad v|^2 - \int_\Omega (m_1 - u + a v) \\
	&\leq& -\frac{\delta}{\us} \int_\Omega \frac{|\grad u|^2}{u^2} + \frac{\chi^2}{4  (d_1 - \frac{\delta}{\us})} \int_\Omega |\grad v|^2 - \int_\Omega (m_1 - u + a v) 
\eas
for all $t\in(0,\infty)$ by integration by parts, an application of Young's inequality as well as the estimate $F(y) \leq y$ for all $y \geq 0$. Combining the above two inequalities then allows us to estimate
\bea{ddt_Hu_ineq} \hspace{-7mm}
	\frac{d}{dt} \int_\Omega H_\us(u) \nn
	&=& \frac{d}{dt} \left( \int_\Omega u - \us - \us \log\left(\frac{u}{\us}\right) \right)\nn \\
	&\leq& -\delta\int_\Omega \frac{|\grad u|^2}{u^2} + \frac{\chi^2 \us}{4 (d_1 - \frac{\delta}{\us})}  \int_\Omega |\grad v|^2 + \int_\Omega (u - \us) (m_1 - u + a v)\nn \\
	&=& -\delta\int_\Omega \frac{|\grad u|^2}{u^2} + \frac{\chi^2 \us}{4 (d_1 - \frac{\delta}{\us})}  \int_\Omega |\grad v|^2 - \int_\Omega (u - \us)^2 + a\int_\Omega (u - \us)(v - \vs) 
\eea
for all $t\in(0,\infty)$ as $(m_1 - \us + a\vs) = 0$.
If $\vs = 0$, integrating the second equation in (\ref{prob-general}) over $\Omega$ similarly yields
\bea{ddt_Hv_ineq}
	\frac{d}{dt} \int_\Omega H_\vs(v) &=&  \frac{d}{dt} \int_\Omega v = \int_\Omega v(m_2 - bu - v) \nn \\
	&=& -\int_\Omega v(v-\vs) - b\int_\Omega v(u - \us) + \int_\Omega v(m_2 - b\us - \vs) \nn \\
	&\leq& -\int_\Omega v(v-\vs) - b\int_\Omega v(u - \us) \nn \\
	&=& -d_2 \vs \int_\Omega \frac{|\grad v|^2}{v^2} -\int_\Omega (v-\vs)^2 - b\int_\Omega (v-\vs)(u - \us) 
\eea
for all $t\in(0,\infty)$ as $(m_2 - b\us - \vs) \leq 0$. The same inequality can be achieved for the case $\vs > 0$ by essentially the same procedure as the one applied to $u$ to gain (\ref{ddt_Hu_ineq}) because in this case $(m_2 - b\us - \vs) = 0$.
\abs
Combining (\ref{ddt_Hu_ineq}) with an appropriately scaled version of (\ref{ddt_Hv_ineq}) then results in the estimate
\bas
	& & \hspace{-13mm} \frac{d}{dt} \left[ \int_\Omega H_\us(u) + \frac{a}{b}\int_\Omega H_\vs(v) \right] \\
	&\leq&  -\delta \int_\Omega \frac{|\grad u|^2}{u^2} +  \frac{\chi^2 \us}{4 (d_1 - \frac{\delta}{\us})} \int_\Omega |\grad v|^2 - d_2 \vs \frac{a}{b}  \int_\Omega \frac{|\grad v|^2}{v^2} - \int_\Omega (u - \us)^2 - \frac{a}{b}\int_\Omega (v - \vs)^2
\eas
for all $t \in (0, \infty)$. By now further testing the second equation in (\ref{prob-general}) with $(v-\vs)$, we gain
\bas
	& & \hspace{-13mm} \frac{1}{2}\frac{d}{dt} \int_\Omega (v - \vs)^2 \\
	&=& -d_2\int_\Omega |\grad v|^2 + \int_\Omega (v - \vs)v(m_2 - bu - v) \\
	&=& -d_2\int_\Omega |\grad v|^2 -\int_\Omega (v-\vs)^2 v - b\int_\Omega (v-\vs)v(u-\us) + \int_\Omega (v - \vs)v(m_2 - b\us - \vs) \\
	&\leq& -d_2\int_\Omega |\grad v|^2 -\int_\Omega (v-\vs)^2 v - b\int_\Omega (v-\vs)v(u-\us) \\
	&\leq& -d_2\int_\Omega |\grad v|^2 + \frac{b^2}{4} \|v\|_{L^\infty(\Omega)} \int_\Omega (u-\us)^2
\eas
for all $t \in (0, \infty)$ because $(m_2 - b\us - \vs) \leq 0$ for $\vs = 0$ and $(m_2 - b\us - \vs) = 0$ for $\vs > 0$. Therefore, combining this with the inequality directly preceding it yields
\bea{stabilization_energy_1}
	\frac{d}{dt}\int_\Omega E(u,v)
	&\leq&  -\delta\int_\Omega \frac{|\grad u|^2}{u^2} + \left(\frac{\chi^2 \us}{4 (d_1 - \frac{\delta}{\us})} - \frac{4d_2}{b^2 \widehat{m_2}}  \right) \int_\Omega |\grad v|^2 -  d_2 \vs \frac{a}{b} \int_\Omega \frac{|\grad v|^2}{v^2} \nn \\
	& & + \left(\frac{\|v\|_{L^{\infty}(\Omega)}}{\widehat{m_2}} - 1\right)\int_\Omega (u - \us)^2 - \frac{a}{b}\int_\Omega (v - \vs)^2 
\eea
for all $t \in (0, \infty)$.
\abs
Notably, the above implies that
\bas
	\frac{d}{dt}\int_\Omega E(u,v) &\leq& \frac{\chi^2 \us}{4 (d_1 - \frac{\delta}{\us})} \int_\Omega |\grad v|^2 + \frac{\|v\|_{L^\infty(\Omega)}}{\widehat{m_2}}\int_\Omega (u - \us)^2
\eas
for all $t \in (0,\infty)$. By application of Lemma \ref{lem_basic-v} and Lemma \ref{lem_basic-u} as well as the initial data condition (\ref{stabilization_initial_data_condition}), which ensures that $E(u_0, v_0) < \infty$, this directly gives us uniform boundedness of the energy functional on $[0, T_E]$ and thus the first half of our result by time integration.  
\abs
If we now consider the inequality (\ref{stabilization_energy_1}) only on the time interval $(T_E, \infty)$, we can use (\ref{energy_eventual_m_bound}) followed by (\ref{energy_param_1})--(\ref{energy_param_3}) to further estimate that
\bas
	\frac{d}{dt}\int_\Omega E(u,v)
	&\leq&  -\delta\int_\Omega \frac{|\grad u|^2}{u^2} + \left[ \left(\frac{\chi^2 \us}{4 (d_1 - \frac{\delta}{\us})} - \frac{4d_2}{b^2 \widehat{m_2}}  \right) \widehat{m_2}^2   -d_2 \vs \frac{a}{b} \right]\int_\Omega \frac{|\grad v|^2}{v^2}  \\
	& & + \left(\frac{\|v\|_{L^{\infty}(\Omega)}}{\widehat{m_2}} - 1\right)\int_\Omega (u - \us)^2 - \frac{a}{b}\int_\Omega (v - \vs)^2  \\
	&\leq& -\delta\int_\Omega \frac{|\grad u|^2}{u^2} - \delta \int_\Omega \frac{|\grad v|^2}{v^2} - \delta \int_\Omega (u - \us)^2  - \delta\int_\Omega (v - \vs)^2 = -\delta G(u,v)
\eas
for all $t \in (T_E, \infty)$, which completes the proof.
\qed
%
%
%
%
%
%
\subsection{Monotonicity of the energy-type functional}
Having now established our energy inequality, the first key consequence we will derive from it is that after the time $T_E$ our energy functional is in fact always monotonically decreasing. We do this to later gain an eventual smallness property for $E$ by just constructing an increasing sequence of times along which $E$ becomes small, which combined with monotonicity is of course sufficient.
\abs
While the aforementioned monotonicity property is a trivial consequence of Lemma \ref{lem_stab_energy} in the two-dimensional case, in the three-dimensional case we still need to argue that the same monoticity property survives the limit process used in our weak solution construction. In fact, given that the constructed weak solutions are merely nonnegative, it is not even immediately clear that $E(u,v)$ is finite for almost all times $t > T_E$ if $(u,v)$ is such a weak solution. As such, we will now first prove what is effectively a uniform lower bound for $\int_\Omega \log(u)$ and $\int_\Omega \log(v)$.
\begin{lem}\label{lem_ln_uv_finite}
	Under the assumption of Lemma \ref{lem_stab_energy}, there exists $C > 0$ such that
	\bas
		\us \int_\Omega |\log(u(\cdot, t))| \leq C \;\;\;\; \text{ and }  \;\;\;\; \vs \int_\Omega |\log(v(\cdot, t))| \leq C
	\eas
	for all $t \in (T_E, \infty)$ with $T_E$ as in Lemma \ref{lem_stab_energy}.
\end{lem}
\proof%
We fix $T_E > 0$, $C > 0$ as provided by Lemma \ref{lem_stab_energy}.
Then using the energy inequality from the same lemma, we can conclude that
\bas
	& &\hspace*{-11mm} -\us \int_\Omega \log(u(\cdot, t)) \\
	&=& -\us \int_\Omega \log\left(\frac{u(\cdot, t)}{\us}\right) - \us \int_\Omega \log(\us) - E(u(\cdot,t),v(\cdot,t)) + E(u(\cdot,t),v(\cdot,t)) \\
	&=&  - \us \int_\Omega \log(\us) - \int_\Omega (u(\cdot,t) - \us) - \frac{a}{b}\int_\Omega H_\vs(v(\cdot,t)) - \frac{2}{b^2 \widehat{m_2}} \int_\Omega (v(\cdot,t) - \vs)^2 + E(u(\cdot, t), v(\cdot,t)) \\
	&\leq& - \us \int_\Omega \log(\us) + \int_\Omega  \us + E(u(\cdot, T_E), v(\cdot,T_E)) \\
	&\leq&  \us (1 + |\log(\us)|)|\Omega|  + E(u(\cdot, T_E), v(\cdot,T_E)) \leq \us (1 + |\log(\us)|)|\Omega| + C
\eas
for all $t \in (T_E, \infty)$. As further
\bas
	\us\int_\Omega |\log(u(\cdot, t))| &=& -\us\int_{\{u \leq 1\}} \log(u(\cdot, t)) + \us\int_{\{u > 1\}} \log(u(\cdot, t)) \\
	&\leq& -\us\int_\Omega \log(u(\cdot, t)) + 2\us\int_{\{u > 1\}} u(\cdot, t)
\eas
for all $t \in (T_E, \infty)$, our first bound follows directly due to Lemma \ref{lem_basic-u}.
The remaining bound for $v$ can be achieved in essentially the same way if $\vs > 0$ and is trivial if $\vs = 0$.
\qed
We further derive a uniform square integrability property for $H_\us(u)$ and $H_\vs(v)$ based on the above result as well as Lemma \ref{lem_stab_energy}. 
\begin{lem}\label{lem_H2_bounds}
	Under the assumption of Lemma \ref{lem_stab_energy}, there exists $C > 0$ such that
	\bas
	\int_t^{t+1}\int_\Omega H_\us^2(u) \leq C \;\;\;\; \text{ and } \;\;\;\; \int_t^{t+1}\int_\Omega H_\vs^2(v) \leq C
	\eas
	for all $t \in (T_E, \infty)$ with $T_E$ as in Lemma \ref{lem_stab_energy}.
\end{lem}
\proof%
Using Lemma \ref{lem_stab_energy} as well as Lemma \ref{lem_ln_uv_finite}, we can fix $T_E > 0$, $c_1 > 0$ and $\delta > 0$ such that
\bas
	& & \hspace*{-13mm} \int_{T_E}^\infty \int_\Omega \frac{|\grad u|^2}{u^2} + \int_{T_E}^\infty \int_\Omega \frac{|\grad v|^2}{v^2} + \int_{T_E}^\infty\int_\Omega (u -\us)^2  + \int_{T_E}^\infty\int_\Omega (v -\vs)^2 \\
	& = &\int_{T_E}^\infty G(u,v) \leq \tfrac{1}{\delta}E(u(\cdot, T_E), v(\cdot, T_E)) \leq \tfrac{c_1}{\delta}
\eas
as well as
\bas
	\us \int_\Omega |\log(u(\cdot, t))| \leq c_1 \;\;\;\; \text{ and } \;\;\;\; \vs \int_\Omega |\log(v(\cdot, t))| \leq c_1
\eas
for all $t \in (T_E, \infty)$. An application of (\ref{Hxi_upper_bound}) from Lemma \ref{lem_H_props} combined with the above bounds then yields $c_2 > 0$ such that 
\bas
	& & \hspace*{-13mm} \int_t^{t+1} \int_\Omega H^2_\us(u) \\
	&\leq& \int_t^{t+1}\int_\Omega u^2 + c_2\left( u_*^2 \int_t^{t+1}\int_\Omega \frac{|\grad u|^2}{u^2} + u_*^2\int_t^{t+1}\left\{ \int_\Omega |\log(u)| \right\}^2 + u_*^2 \right) \\
	&\leq& 2\int_t^{t+1}\int_\Omega (u - \us)^2 + 2|\Omega|u_*^2 + c_2\left( u_*^2 \int_t^{t+1}\int_\Omega \frac{|\grad u|^2}{u^2} + u_*^2\int_t^{t+1}\left\{ \int_\Omega |\log(u)| \right\}^2 + u_*^2 \right) \\
	&\leq& 2 \tfrac{c_1}{\delta} + 2|\Omega| u_*^2 + c_2 \left( u_*^2\tfrac{c_1}{\delta} + c_1^2 + u_*^2\right) 
\eas
for all $t \in (T_E, \infty)$.
Naturally, the same estimation applies to $\int_t^{t+1} \int_\Omega H^2_\vs(v)$, completing the proof.
\qed
Given this, we can now derive our desired (almost everywhere) monotonicity property for $E(u,v)$ in both two and three dimensions.
\begin{lem}\label{lem_E_monotone}
	Let $(u,v)$ be the solution constructed in Theorem \ref{theo_global}. If (\ref{stabilization_initial_data_condition}) as well as (\ref{stablization_param_condition}) hold, then there exists a null set $N\subset (T_E, \infty)$ such that
	\bas
		E(u(\cdot, t_1),v(\cdot, t_1)) \leq E(u(\cdot, t_0), v(\cdot, t_0))
	\eas
	for all $t_0, t_1 \in (T_E, \infty) \setminus N$ with $t_1 > t_0$ and $T_E$ as in Lemma \ref{lem_stab_energy}. If $n = 2$, then $N = \emptyset$.
\end{lem}
\proof%
	If $n = 2$, this is a trivial consequence of Lemma \ref{lem_stab_energy}. As such, we will focus our efforts here on the case $n = 3$ where Lemma \ref{lem_stab_energy} only implies that there exists a uniform $T_E > 0$ such that
	\be{eq:stab_approx_energy_monotone}
		E(\ueps(\cdot, t_1),\veps(\cdot, t_1)) \leq E(\ueps(\cdot, t_0), \veps(\cdot, t_0))
	\ee
	for all $t_1 > t_0 > T_E$ and $\eps \in (0,1)$ with $(\ueps, \veps)$ being the approximate solutions used in Section \ref{sec: 3D}. Therefore, we still need to argue that this property in fact translates to their limit functions for it to apply to the weak solutions constructed in the same section. We first note that due to the convergence properties already laid out in Lemma \ref{lem_limit-process}, we can ensure that there exists a null sequence $(\eps_j)_{j\in\N}$ such that
	\be{eq:stab_ve_L2_convergence}
		\veps(\cdot, t) \rightarrow v(\cdot, t) \;\;\;\; \text{ in } L ^2(\Omega) \;\;\;\; \text{ for a.e.} \;\;\;\;  t \in (T_E, \infty)
	\ee
	and thus
	\bas
	\int_\Omega (\veps(\cdot, t) - \vs)^2 \rightarrow \int_\Omega (v(\cdot, t) - \vs)^2 \;\;\;\; \text{ for a.e.} \;\;\;\;  t \in (T_E, \infty)
	\eas
	as well as $\ueps \rightarrow u$ and $\veps \rightarrow v$ almost everywhere in $\Omega \times (T_E, \infty)$.  
	\abs
	To now treat the remaining parts of the energy functional $E$, we first note that Lemma \ref{lem_ln_uv_finite} gives us $c_1 > 0$ such that
	\bas
		\int_\Omega |\log(\ueps(\cdot, t))| \leq c_1
	\eas
	for all $t \in (T_E,\infty)$ and $\eps \in (0,1)$ as $\us > 0$. Due to Fatou's lemma this implies  
	\bas
		\int_\Omega |\log(u(\cdot, t))| \leq c_1
	\eas
	for all $t \in (T_E,\infty)$ and thus $u$ must be positive for almost all $(x,t) \in \Omega \times (T_E, \infty)$. As such and due to the fact that $H_\us$ is continuous outside of $0$, we know that $H_\us(\ueps) \rightarrow H_\us(u) $ almost everywhere in $\Omega \times (T_E, \infty)$ as $\eps = \eps_j \searrow 0$ due to the already established almost everywhere convergence for $\ueps$ itself. Given that Lemma \ref{lem_H2_bounds} grants us a constant $c_2 > 0$ with
	\[
		\int_t^{t+1}\int_\Omega H_\us^2(\ueps) \leq c_2
	\]
	for all $t \in (T_E,\infty)$ and $\eps \in (0,1)$, the Vitali convergence theorem gives us that further
	\[
		\int_t^{t+1}\int_\Omega H_\us(\ueps) \rightarrow \int_t^{t+1}\int_\Omega H_\us(u)
	\] as $\eps = \eps_j \searrow 0$ for all $t \in (T_E,\infty)$. But this directly implies 
	\[
		\int_\Omega H_\us(\ueps(\cdot, t)) \rightarrow \int_\Omega H_\us(u(\cdot, t))
	\]	
	for almost all $t \in (T_E,\infty)$. If $\vs > 0$, the argument for $\int_\Omega H_\vs(v)$ is essentially the same while for $\vs = 0$ it is a straightforward consequence of (\ref{eq:stab_ve_L2_convergence}). As such, the property (\ref{eq:stab_approx_energy_monotone}) translates to the limit functions outside of a set $N \subset (T_E, \infty)$ with measure zero by simply taking the limit $\eps = \eps_j \searrow 0$. 
\qed

%
%
%
%
%
%
\subsection{Eventual smallness of the energy functional and proof of Theorem \ref{theo_stabilization}}
As already alluded to in the previous section, we will now construct a sequence of times along which $E(u,v)$ becomes small. To do this, we now first translate some important space-time integral bounds stemming from Lemma \ref{lem_stab_energy} and Lemma \ref{lem_H2_bounds} to the solutions constructed in Section \ref{sec: global} using Fatou's lemma.
\begin{lem}\label{lem_stab_conv_integral_bounds}
	Let $(u,v)$ be the solution constructed in Theorem \ref{theo_global}. If (\ref{stabilization_initial_data_condition}) as well as (\ref{stablization_param_condition}) hold, then there exists $C > 0$ such that 
	\bas
		\int_{T_E}^\infty \int_\Omega (u - \us)^2 \leq C \;\;\;\; \text{ and } \;\;\;\; \int_{T_E}^\infty \int_\Omega (v - \vs)^2 \leq C 
	\eas
	as well as
	\bas
		\int_t^{t+1}\int_\Omega H^2_\us(u) \leq C \;\;\;\; \text{ and } \;\;\;\; \int_t^{t+1} \int_\Omega H^2_\vs(v) \leq C
	\eas
	for all $t \in (T_E, \infty)$ with $T_E$ as in Lemma \ref{lem_stab_energy}.
\end{lem}
\proof%
	If $n = 2$, this is a trivial consequence of Lemma \ref{lem_stab_energy} and Lemma \ref{lem_H2_bounds}. If $n = 3$, Fatou's lemma combined with the same lemmata yields the desired result.
\qed
Due to the above combined with the Vitali convergence theorem, it now follows that the energy $E$ becomes small along a sequence of increasing times, which combined with the monotonicity property of the previous section and Lemma \ref{lem_H_props} leads to the following stabilization properties.
\begin{lem}\label{lem_stabilization}
	Let $(u,v)$ be the solution constructed in Theorem \ref{theo_global}.  If (\ref{stabilization_initial_data_condition}) as well as (\ref{stablization_param_condition}) hold, then there exists a null set $N \subset (0,\infty)$ such that
	\bas
		u(\cdot, t) \rightarrow \us \;\;\;\; \text{ in } \;\;\;\; L^1(\Omega)
	\eas
	and
	\bas
		v(\cdot, t) \rightarrow \vs \;\;\;\; \text{ in } \;\;\;\; L^p(\Omega) \;\;\;\; \text{ for all } p \in [1,\infty)
	\eas
	as $(0,\infty)\setminus N \ni t \rightarrow \infty$. If $n = 2$, then $N = \emptyset$.
\end{lem}
\proof%
	As a direct consequence of Lemma \ref{lem_stab_conv_integral_bounds}, we can fix a sequence of times $(t_k)_{k\in\N} \subset (T_E,\infty)\setminus N$ with $t_k \nearrow \infty$ as $k\rightarrow \infty$, $T_E$ as in Lemma \ref{lem_stab_energy} and $N$ as in Lemma \ref{lem_E_monotone} as well as a constant $c > 0$ such that
	\be{stab_subseq_l2_conv}
		\int_\Omega( u(\cdot, t_k) - \us )^2 \rightarrow 0 \;\;\;\; \text{ and } \;\;\;\; \int_\Omega( v(\cdot, t_k) - \vs )^2 \rightarrow 0 
	\ee
	as $k\rightarrow \infty$ and
	\be{stab_subseq_H2_bound}
		\int_\Omega H^2_\us(u(\cdot, t_k)) \leq c \;\;\;\; \text{ and } \;\;\;\;\int_\Omega H^2_\vs(v(\cdot, t_k)) \leq c
	\ee
	for all $k \in \N$. By potentially switching to a subsequence, we can further ensure
	that $u(\cdot, t_k) \rightarrow \us$ and $v(\cdot, t_k) \rightarrow \vs$ almost everywhere as $k \rightarrow \infty$ as a direct consequence of the $L^2(\Omega)$ convergence established in (\ref{stab_subseq_l2_conv}). As notably $H_\xi$ is continuous in $\xi$ for all $\xi \geq 0$, this further implies that
	\bas
		H_\us(u(\cdot, t_k)) \rightarrow H_\us(\us) = 0 \;\;\;\; \text{ and } \;\;\;\; H_\vs(v(\cdot, t_k)) \rightarrow H_\vs(\vs) = 0
	\eas
	almost everywhere as $k \rightarrow \infty$. Because we have already established the bound seen in (\ref{stab_subseq_H2_bound}), we can now apply the Vitali convergence theorem to conclude that
	\bas
		\int_\Omega H_\us(u(\cdot, t_k)) \rightarrow 0 \;\;\;\; \text{ and } \;\;\;\; \int_\Omega H_\vs(v(\cdot, t_k)) \rightarrow 0
	\eas
	as $k \rightarrow \infty$ as well. Combined with (\ref{stab_subseq_l2_conv}), this implies that
	\bas
		E(u(\cdot, t_k), v(\cdot, t_k)) \rightarrow 0
	\eas
	as $k \rightarrow \infty$, which due to the monotonicity property laid out in Lemma \ref{lem_E_monotone} gives us that 
	\bas
		E(u(\cdot, t), v(\cdot, t)) \rightarrow 0
	\eas
	as $(0, \infty)\setminus N \ni t \rightarrow \infty$. Then an application of inequality (\ref{Hxi_lower_bound}) from Lemma \ref{lem_H_props} yields that $u(\cdot, t) \rightarrow \us$ and $v(\cdot, t) \rightarrow \vs$ in $L^1(\Omega)$ as $(0, \infty)\setminus N \ni t \rightarrow \infty$, which combined with the fact that $v$ is globally bounded in $L^\infty(\Omega)$ due to Theorem \ref{theo_global} completes the proof.
\qed%
To complete this section, we will now use the regularity properties already established in Theorem \ref{theo_global} combined with a brief testing based argument to improve upon the strength of the above stabilization properties in the two-dimensional case and thus prove our second central theorem.
\abs
\proofc of Theorem \ref{theo_stabilization}.\quad
For $n = 3$, Theorem \ref{theo_stabilization} is an immediate consequence of Lemma \ref{lem_stabilization}. 
\abs
As such, we will now focus on the case $n = 2$. Here, we first note that as consequence of the properties derived in Lemma \ref{lem_stabilization} we gain
\bas
	u(\cdot, t) \rightarrow \us  \;\;\;\; \text{ in } \;\;\;\; L^p(\Omega)\;\;\;\; \text{ for all } p \in [1,\infty)	
\eas
as $t \rightarrow \infty$ due to the $L^\infty(\Omega)$ bound for $u$ established in Theorem \ref{theo_global} if $n = 2$. Similar to Lemma \ref{lem_Lap-v}, we now test the second equation in (\ref{prob}) with $-\Delta v$ and employ integration by parts as well as Young's inequality to gain
\bas
	\frac{1}{2}\frac{d}{dt} \int_\Omega |\grad v|^2 + d_2 \int_\Omega |\Delta v|^2 &=& -\int_\Omega v(m_2-bu-v) \Delta v \\
	&=& b\int_\Omega v(u - \us)\Delta v +\int_\Omega v (v - \vs) \Delta v - \int_\Omega v (m_2 - b\us - \vs) \Delta v \\
	&=& b\int_\Omega v(u - \us)\Delta v +\int_\Omega v (v - \vs) \Delta v + \int_\Omega (m_2 - b\us - \vs) |\grad v|^2 \\
	&\leq& \frac{d_2}{2} \int_\Omega |\Delta v|^2 + \frac{b^2}{d_2}\|v\|^2_{L^{\infty}(\Omega))} \int_\Omega (u - \us)^2 + \frac{1}{d_2}\|v\|^2_{L^{\infty}(\Omega))} \int_\Omega (v - \vs)^2
\eas
for all $t \in (0,\infty)$ as $(m_2 - b\us - \vs) \leq 0$ . Given that we have already established that $u(\cdot, t)$ and $v(\cdot, t)$ converge to $\us$ and $\vs$ in $L^2(\Omega)$ as $t\rightarrow \infty$ and $v$ is globally bounded in $L^\infty(\Omega)$ as well as by application of a well-known consequence of the Poincar\'e inequality when considering functions with Neumann boundary conditions, it immediately follows that there exists a constant $c > 0$ and, for each $\delta > 0$, a time $t_0 = t_0(\delta) > T_E$ with $T_E$ as in Lemma \ref{lem_stab_energy} such that
\bas
	\frac{d}{dt} \int_\Omega |\grad v|^2 + c \int_\Omega |\grad v|^2 \leq \frac{\delta}{2} c
\eas
for all $t \in (t_0, \infty)$. By a standard comparison argument this implies that 
\bas
	\int_\Omega |\grad v(\cdot, t)|^2 \leq \left(\int_\Omega |\grad v(\cdot, t_0)|^2 - \frac{\delta}{2}\right)e^{-c (t-t_0)} + \frac{\delta}{2}
\eas
for all $t \in (t_0, \infty)$ and therefore after an appropriate waiting time $t_1 = t_1(\delta) > t_0$ it follows that
\[
	\int_\Omega |\grad v(\cdot, t)|^2 \leq \delta
\]
for all $t \in (t_1, \infty)$. As such, we have gained that
\bas
	v(\cdot, t) \rightarrow \vs  \;\;\;\; \text{ in } \;\;\;\; W^{1,2}(\Omega)	
\eas
as $t \rightarrow \infty$, which combined with the $W^{1,\infty}(\Omega)$ bound for $v$ from Theorem \ref{theo_global}, completes the proof.
\qed

\bigskip

{\bf Acknowledgment.}
This work was started while the second author was supported by 
Grant-in-Aid for Scientific Research (C), No.\,16K05182, JSPS, 
and carried out while he visited Universit\"at Paderborn in 2019-2020 by 
Tokyo University of Science Grant for International Joint Research. The first author further acknowledges the support of the Deutsche Forschungsgemeinschaft within the scope of the
project Emergence of structures and advantages in cross-diffusion systems, project number 411007140.

\end{document}